\documentclass[11pt,a4paper,twoside,english,notitlepage]{article}
\usepackage{amsmath}
\usepackage{amssymb}

\makeatletter
\usepackage{graphicx}
\usepackage{amsfonts}
\newtheorem{theorem}{Theorem}

\newtheorem{corollary}[theorem]{Corollary}

\newtheorem{lemma}[theorem]{Lemma}

\newtheorem{proposition}[theorem]{Proposition}
\newtheorem{remark}[theorem]{Remark}

\numberwithin{theorem}{section}
\newenvironment{proof}[1][Proof]{\textbf{#1.} }{\ \rule{0.5em}{0.5em}}
\DeclareMathOperator{\Ad}{Ad}
\DeclareMathOperator{\ad}{ad}

\usepackage{babel}
\makeatother
\begin{document}

\title{A Cotangent Bundle Slice Theorem
%\\
%\textsc{working draft- not for circulation}
}

\author{Tanya Schmah
\footnote{
Department of Mathematics,
Macquarie University,
NSW 2122 Australia,
\texttt{schmah@maths.mq.edu.au},
\texttt{http://www.maths.mq.edu.au/\~{}schmah}
}}

\date{20 August 2004}

\maketitle
\begin{abstract}
This article concerns cotangent-lifted Lie group actions;
our goal is to find local and ``semi-global'' normal forms for these and associated structures. 
Our main result is a constructive cotangent bundle slice
theorem that extends the Hamiltonian slice theorem of Marle \cite{Marl85}
and Guillemin and Sternberg \cite{GS84}.
The result applies to all proper cotangent-lifted actions, around points with
fully-isotropic momentum values.

We also present a ``tangent-level'' commuting reduction result and use it to 
characterise the
symplectic normal space of any cotangent-lifted action.
In two special cases, we arrive at splittings of the symplectic normal space.
One of these cases is when the configuration isotropy group is contained in the
momentum isotropy group; in this case, our splitting generalises that 
given for free actions by Montgomery et al. \cite{MMR84}.
The other case includes all relative equilibria of simple mechanical systems. 
In both of these special cases, the new splitting leads to a refinement of the so-called \emph{reconstruction equations} or \emph{bundle equations} \cite{Ort98,OR02b,RWL02}.
We also note cotangent-bundle-specific local normal forms for symplectic reduced spaces.
\end{abstract}
\tableofcontents{}

\section{Introduction}

This article concerns cotangent-lifted actions of a Lie group $G$ on a cotangent bundle $T^*Q.$
We are motivated in part by the role of such actions as groups of symmetries of
Hamiltonian systems with cotangent bundle phase spaces.
Nonetheless, this article is primarily geometric, the exception being the discussion of the
reconstruction equations (bundle equations) at the end of Section \ref{sectNs}.

When $G$ acts freely and properly, it is well known that
one can \emph{reduce} $T^*Q$ by the $G$ action to give a lower-dimensional symplectic manifold
(see Theorem \ref{RPSR}). The reduced manifold inherits some cotangent-bundle structure
\cite{AM78,Mars92,MP00}, and it sometimes actually \emph{is} a cotangent bundle 
(see for example Theorem \ref{RPCBRZ}).
In Hamiltonian systems, the solutions of the system on the original space project to 
solutions of a new Hamiltonian system on the reduced phase space.

One would like to generalise this picture to arbitrary proper group actions, not necessarily free.
This problem of \emph{singular reduction} has been addressed with success
in the symplectic category \cite{SL91,BL97} and more recently for the special case of 
cotangent bundles \cite{CS01, PRS03}. But the symplectic reduced spaces are in general
not smooth, and our understanding of the inherited cotangent bundle structure is far from
complete. 

A different but related approach is to ask: to what degree can we factor out the symmetry
while not losing smoothness? Slice theorems are an answer to this question.
For any \emph{free} proper action of $G$ on a manifold $M,$
the slice theorem of Palais (Theorem \ref{slice}) says that every point $z\in M$ has a neighbourhood
$G$-equivariantly isomorphic to a space $G\times S,$ where $S$ is some submanifold of $M$
transverse to the $G$ orbit and the $G$ action on $G\times S$ is $g'\cdot(g,s)=(g'g,s).$
This local model of the action of $G$ on $M$ is actually ``semi-global''  in the sense that it is
global in the $G$ direction but local in the transverse direction.
For general proper actions, the model space is not $G\times S$ but a \emph{twisted} space
$G\times_{G_z} S,$ where $G_z$ is the isotropy group of the point $z.$
For symplectic actions, the Hamiltonian slice theorem of Marle \cite{Marl85}
and Guillemin and Sternberg \cite{GS84} (Theorem \ref{marle})
gives a model space of this kind and a $G$-equivariant \emph{symplectic} diffeomorphism.
This theorem is a fundamental tool in the study of Hamiltonian systems with symmetry:
it has found applications to singular reduction \cite{SL91,BL97,Ort98}, and to many dynamical
questions involving stability, bifurcation and persistence in the
neighbourhood of relative equilibria and relative periodic orbits
\cite{Mo97,Mo97a,RdSD97,LS98,OR99,RWL02,Ort03,OR02a,WR02}. 

The main aim of the present article is to extend the Hamiltonian slice theorem in the
context of cotangent bundles. We succeed in doing so around
points with fully isotropic momentum values (Theorem \ref{mainGmuG}).
Our new result extends
that of Marle, Guillemin and Sternberg in three ways.
First, it involves a new cotangent-bundle-specific splitting of the symplectic
normal space. 
Second, it is constructive, up to a Riemannian exponential map.
In particular, we do not use the Constant Rank Embedding Theorem or
Darboux's Theorem. 
Third, our construction has a unique property (see Lemma \ref{sigma}).

We begin with a summary of the relevant background material, including some reduction
results and slice theorems.
In Section \ref{sectcommred} we summarise regular and singular commuting symplectic reduction,
and introduce a new ``tangent level'' commuting reduction result that works at the level
of symplectic normal spaces (Theorem \ref{TLCR}).

In Section \ref{sectNs} we analyse the symplectic normal space at $z\in T^*_qQ$ 
for a cotangent-lifted action of $G$ on $T^*Q.$
We first apply Palais' slice theorem in the configuration space and then 
cotangent-lift the resulting diffeomorphism to give a local symplectic diffeomorphism from  $T^*Q$ 
to $T^*\left(G\times_{G_q} A\right),$
where $A$ is a $G_z$-invariant complement to the tangent to the orbit $G\cdot q.$
We then ``unroll'' this space by considering the untwisted product $T^*\left(G\times A\right).$
We note that there are two obvious commuting actions on this space, namely 
cotangent lifts of left multiplication by $G$ and twist by $G_q.$
Applying tangent-level commuting reduction leads to our characterisation of the
symplectic normal space $N_s$ of the original $G$ action on $T^*Q,$ Theorem \ref{mainVm}.
Corollaries \ref{NsKsubGmu} and \ref{Nsalpha0} give
splittings of $N_s$ in two special cases: $G_z \subset G_\mu,$ where $\mu$ is the momentum 
value of $z; $ and $z$ ``purely in the group direction'', meaning $\left. z\right|_A = 0.$
We note consequences of these results for singular reduction and for the
reconstruction equations (bundle equations).

In Section \ref{sectcbs}, we prove the cotangent bundle slice theorem.
We begin with the observation that, when $G_\mu = G,$ 
our new splitting $N_s\cong T^*B$ (for a certain subspace $B$ of $A$)
 implies that the model space in the Hamiltonian slice theorem
is $G\times_{G_z} \left(\left(\mathfrak{g}/\mathfrak{g}_z\right)^* \oplus T^*B\right),$
which is a $G_z$-reduced space of $T^*\left(G\times B\right).$
The problem of proving a constructive Hamiltonian slice theorem thus reduced in this
case to that of finding a certain symplectic local diffeomorphism from
$T^*\left(G\times_{G_z} B\right)$ to
$T^*\left(G\times_{G_q} A\right).$
We proceed to find a suitable symplectomorphism, using two alternative methods.
The first method is more ``brute-force'' and gives an explicit formula in coordinates;
the second method is to re-arrange the problem so that a cotangent-lift can be used.
We end with a simple example, $SO(3)$ acting on $T^*\mathbf{R}^3.$

Most of the results in this article first appeared in the author's PhD thesis
\cite{Sch02}, where more detailed proofs of some results appear.

\bigskip
\noindent
\textbf{Some assumptions and notation:} We consider only proper
actions on finite-dimensional manifolds. All cotangent bundles are
given the standard cotangent bundle symplectic form; in particular,
no magnetic terms appear. All group actions are left actions.
The Lie algebra of a Lie group will always
be denoted by the corresponding fraktur letter. 
%The notation $\pi_{G},$
%for a group $G$ acting on a manifold $M,$ will denote the projection
%$\pi_{G}:M\rightarrow M/G,$ or sometimes a restriction of this map;
%the precise definition should be clear from context. 

\section{Preliminaries}
\label{sectprelim}

We summarise relevant basic facts on Lie group symmetries
symplectic reduction and slice theorems.
This material is well-known;
good general
references are \cite{AM78,DK99,CB97,OR04}. 

\bigskip
\noindent \textbf{Lie Group Actions}
Let $G$ be a Lie group, with Lie algebra $\mathfrak{g},$ and consider
a smooth left action $\Phi$ of $G$ on a manifold $M;$ we write
$g\cdot z=\Phi_{g}\left(z\right)=\Phi\left(g,z\right).$ For every
$\xi\in\mathfrak{g,}$ the \emph{infinitesimal generator} of $\xi$
is the vector field $\xi_{M}$ defined by $\xi_{M}\left(z\right)=\frac{d}{dt}\left.\exp\left(t\xi\right)\cdot z\right|_{t=0}.$
We will also write $\xi_{M}\left(z\right)$ as $\xi\cdot z,$ and
refer to the map $\left(\xi,z\right)\mapsto\xi\cdot z$ as the \emph{infinitesimal
action} of $\mathfrak{g}$ on $M.$

The action $\Phi$ is \emph{proper} if the map $\left(g,z\right)\longmapsto\left(z,g\cdot z\right)$
is proper (i.e. the preimage of every compact set is compact). It
is easily shown that all proper actions on manifolds have the following
property (often used as a definition of properness): given any convergent
sequences $\left\{ z_{i}\right\} $ and $\left\{ g_{i}\cdot z_{i}\right\} ,$
the sequence $\left\{ g_{i}\right\} $ has a convergent subsequence. 

If $G$ acts properly and freely on $M,$ then $M/G$ has a unique
smooth structure such that $\pi_{G}:M\rightarrow M/G$ is a submersion
(in fact, $\pi_{G}$ is a principal bundle). One useful consequence
is that for every $z\in M,$ we have $\ker T_{z}\pi_{G}=T_{z}\left(G\cdot z\right)=\mathfrak{g}\cdot z.$

The \emph{isotropy subgroup} of a point $z\in M$ is $G_{z}:=\left\{ g\in G\mid g\cdot z=z\right\} .$
The isotropy subgroups are always Lie subgroups, since they are clearly
closed. 
An action is \emph{free} if all of the isotropy
subgroups $G_{z}$ are trivial. 

A key elementary property of proper actions is that all isotropy
subgroups are compact. This property is often used produce $H$-invariant
structures by averaging over some isotropy group $H.$

\bigskip
\noindent
\textbf{Momentum Maps}
Now suppose $G$ acts symplectically on a symplectic manifold $\left(M,\omega\right).$ 
Recall that
any function $F:M\rightarrow\mathbf{R}$ defines a Hamiltonian vector
field $X_{F}$ by $i_{X_{F}}\omega=dF,$ in other words $\omega\left(X_{F}\left(z\right),v\right)=dF\left(v\right)$
for every $v\in T_{z}^{\ast}M.$ A \emph{momentum map} is a function
$J:M\rightarrow\mathfrak{g}^{\ast}$ satisfying $X_{J_{\xi}}=\xi_{M}$
for every $\xi\in\mathfrak{g,}$ where $J_{\xi}:M\rightarrow\mathbf{R}$
is defined by $J_{\xi}\left(z\right)=\left\langle J\left(z\right),\xi\right\rangle .$
The $\Ad^{\ast}$ action of $G$ on $\mathfrak{g}^{\ast}$
is given by $g\cdot\nu=\Ad_{g^{-1}}^{\ast}v:=\left(\Ad_{g^{-1}}\right)^{\ast}\nu,$
where $\Ad_{g}$ is the adjoint operator. If the
$G$ action has an $\Ad^{\ast}$-equivariant momentum
map $J,$ then it is called \emph{globally Hamiltonian.}
Note that if $J(z)=\mu$ and $G_{\mu}$ is the isotropy group of $\mu$
with respect to the coadjoint action and $J$ is $\Ad^*$-equivariant then $G_z\subset G_{\mu}.$

We note here an important momentum map, that of the coadjoint action of $G$ on
any coadjoint orbit $\mathcal{O}\subset \mathfrak{g}^*.$
The Kostant-Kirillov-Souriau (KKS) symplectic forms on the coadjoint orbit are
\begin{align}
\label{E:KKS}
\omega_{\mathcal{O}}^{\pm}\left(  \nu\right)  \left(  \xi\cdot\nu
,\eta\cdot\nu\right)  =\pm\left\langle \nu,\left[  \xi,\eta\right]
\right\rangle ,
\end{align}
where \(\xi\cdot\nu=-\ad_{\xi}^{\ast}\nu,\) the
infinitesimal generator of the coadjoint action of \(G\) on \(\mathfrak{g}^{\ast}.\)
The momentum map of the coadjoint action of $G$ on $\mathcal{O}$
with respect to $\omega_{\mathcal{O}}^{\pm}$
is
$J_{\mathcal{O}}\left(\nu\right) = \pm \nu.$

\bigskip
\noindent \textbf{Lifted Actions on (Co-)tangent Bundles.}
Every cotangent bundle $T^{\ast}Q$ has a canonical symplectic form, 
given in given local coordinates by
$\omega=\mathrm{d}q^i\wedge\mathrm{d}p_i.$
The space $Q$ is called the \emph{configuration space} or \emph{base space}.
The \emph{tangent lift} of any action $\Phi:G\times Q\to Q$ is the
action of $G$ on $TQ$ given by $g\cdot v=T\Phi_{g}\left(  v\right)  .$
The \emph{cotangent lift} is the action of $G$ on $T^{\ast}Q$ given by
$g\cdot\alpha=\left(T\Phi_{g^{-1}}\right)^*\alpha.$
Here are some key elementary facts about lifted actions: 

\begin{enumerate}
\item The tangent or cotangent lift of a proper (resp. free) action is proper
(resp. free). 
\item Every cotangent-lifted action on $T^{\ast}Q$ is symplectic with respect
to the canonical symplectic form on $T^{\ast}Q.$
\item \label{defJcb}Every cotangent-lifted action has an $\mathrm{Ad}^{\ast}$-equivariant
momentum map given by 
\[
\left\langle J\left(\alpha_{q}\right),\xi\right\rangle =\left\langle \alpha_{q},\xi\cdot q\right\rangle .
\]
(When we refer to {}``the'' momentum map for such an action, this is the one we mean.) 
%\item If $G$ acts by cotangent lifts on $T^{\ast}Q,$ and if $z\in T_{q}^{\ast}Q,$
%then the isotropy groups of $z$ and $q$ satisfy $G_{z}\subset G_{q},$
%but they are in general not equal. The same is true of tangent lifts. 
%\item If $G$ acts by cotangent lifts on $T^{\ast}Q,$ and if $z\in T_{q}^{\ast}Q$
%and $\mu=J\left(z\right)$ (with $J$ as in point \ref{defJcb}),
%then we have $G_{z}\subset G_{\mu},$ but in general $G_{z}$ may
%be a proper subset of $G_{q}\cap G_{\mu}.$
\end{enumerate}
%The interaction between the different isotropy subgroups $G_{q},G_{z}$
%and $G_{\mu}$ plays a large role in the theory of cotangent-lifted
%actions. 

We now consider the special case where $Q$ is a vector space.
Any $G$ action on a vector space $V$ induces an \emph{inverse dual}
(or \emph{contragredient}) action of $G$ on $V^{\ast}.$ Specifically,
if $\Psi:G\times V\rightarrow V$ is the $G$ action on $V,$ then
the $G$ action on $V^{\ast}$ is $g\cdot\alpha=\left(\Psi_{g^{-1}}\right)^{\ast}\cdot\alpha$.
Identifying $T^{\ast}V$ \ with $V\oplus V^{\ast},$ the cotangent-
lifted action of $V$ on $T^{\ast}V$ is $g\cdot\left(a,\alpha\right)=\left(g\cdot a,g\cdot\alpha\right),$
where the action on the second component is the inverse dual action.
The infinitesimal action of $G$ on $V^{\ast}$ is $\left\langle \eta\cdot\alpha,a\right\rangle =\left\langle \alpha,-\eta\cdot a\right\rangle .$

Note that if we identify $V^{**}$ with $V$ then the inverse dual of the inverse dual of an action
is the original action.

We introduce the diamond notation of Holm et al. \cite{HMR98}, adding an optional subscript to specify
the Lie algebra of the symmetry group or some linear subspace of it. 
For every \(a\in V,\) \(\alpha\in V^{\ast}\), 
and any subspace \(\mathfrak{l}\) of
\(\mathfrak{g}\), we define
\(a\diamond\alpha\in\mathfrak{g}^{\ast}\) and \(a\diamond_\mathfrak{l}\alpha\in\mathfrak{l}^{\ast}\) by 
\begin{align}\label{E:diamond}
\left\langle a\diamond \alpha,\xi\right\rangle &=\left\langle \alpha,\xi\cdot
a\right\rangle 
\textrm{ for all } \xi\in\mathfrak{g}, \textrm{ and }  \\
a\diamond_{\mathfrak{l}}\alpha &=
\left.a\diamond \alpha\right|_{\mathfrak{l}} \notag
\end{align}
We will have occasion
to use the isomorphism $T^{\ast}V\cong V\oplus V^{\ast}\cong V^{\ast\ast}\oplus V^{\ast}\cong T^{\ast}V^{\ast},$
so we point out that $\alpha\diamond_{\mathfrak{g}}a=-a\diamond_{\mathfrak{g}}\alpha$
for any $a\in V$ and $\alpha\in V^*.$

The momentum map for the cotangent-lifted $G$ action on $T^{\ast}V\cong V\oplus V^{\ast}$
is the map $J_{V}:V\rightarrow\mathfrak{g}^{\ast}$ given by $J_{V}\left(a,\alpha\right)=a\diamond_{\mathfrak{g}}\alpha.$

\bigskip
\noindent \textbf{Symplectic Reduction}
We now present three symplectic reduction theorems.
The simplest one is linear:

\begin{theorem}\label{linred}
Let $\omega$ be a skew-symmetric bilinear form on a vector space $V.$
Then $V/V^{\omega}$ has a symplectic form given by
$\omega_{red}\left(u+V^{\omega},v+V^{\omega}\right)=\omega\left(u,v\right).$
\end{theorem}

For symplectic manifolds, one of the simplest forms of reduction is the following,
which is a slight simplification of a version due to Marsden and Weinstein \cite{MW74}.

\begin{theorem}
[Regular ``Point'' Symplectic Reduction]\label{RPSR}Let \(G\) act
freely, properly and symplectically on a symplectic manifold \(\left(
M,\omega\right)  ,\) with \(Ad^{\ast}\)-equivariant momentum map \(J:M\rightarrow
\mathfrak{g}^{\ast},\) and let \(\mu\) \(\in\mathfrak{g}^{\ast}.\) Then the
\emph{reduced space} \(J^{-1}\left(  \mu\right)  /G_{\mu}\) has a symplectic
form \(\omega_{\mu}\) uniquely defined by \(\pi_{\mu}^{\ast}\omega_{\mu}=i_{\mu
}^{\ast}\omega,\) where \(\pi_{\mu}:J^{-1}\left(  \mu\right)  \rightarrow
J^{-1}\left(  \mu\right)  /G_{\mu}\) \(\ \)and \(i_{\mu}:J^{-1}\left(  \mu\right)
\rightarrow M\) is inclusion.
\end{theorem}

%\begin{remark}
%In place of \(G\) acting freely and properly, the following weaker conditions
%suffice: \(\mu\) a regular value of \(J\) and \(G_{\mu}\) acting freely and properly
%on \(J^{-1}\left(  \mu\right)  .\) 
%We will not need this stronger version.
%\end{remark}

If the symplectic manifold is a cotangent bundle, then the reduced
space takes a special form. The following version of cotangent bundle reduction is 
a special case of results by Satzer and Marsden (see \cite{AM78}).

\begin{theorem}
[Regular ``Point'' Cotangent Bundle Reduction at Zero]\label{RPCBRZ}
Let \(G\)
act freely and properly by cotangent lifts on \(T^{\ast}Q,\) and let \(J\) be the
momentum map of the \(G\) action (with respect to the canonical symplectic form
on \(T^{\ast}Q).\) Then \(J^{-1}\left(  0\right)  \) is a smooth submanifold of
\(Q.\) Let \(\pi_{G}:Q\rightarrow Q/G\) be projection. Define the map
\(\varphi:J^{-1}\left(  0\right)  \rightarrow T^{\ast}\left(  Q/G\right)  \) by,
for every \(q\in Q\) and every \(p\in T_{q}^{\ast}Q\) and \(v\in T_{q}Q,\)
\[
\left\langle \varphi\left(  p\right)  ,T\pi_{G}\left(  v\right)  \right\rangle
=\left\langle p,v\right\rangle .
\]
Then \(\varphi\) is a \(G\)-invariant surjective submersion and drops to a
symplectomorphism (i.e. symplectic diffeomorphism)
\[
\bar{\varphi}:J^{-1}\left(  0\right)  /G\rightarrow T^{\ast}(Q/G),
\]
where the left-hand side has the reduced symplectic form corresponding to the
canonical symplectic form on \(T^{\ast}Q,\) and \(T^{\ast}(Q/G)\) has the standard
symplectic form.
\end{theorem}

%\begin{remark}
%The reader familiar with cotangent bundle reduction may recall that in general
%the symplectic form on \(T^{\ast}(Q/G)\) is the canonical one plus a new
%``magnetic'' term; but the latter is always zero at zero momentum.
%\end{remark}

%\begin{remark}
The map \(\varphi\) is a sort of push-forward, though \(\pi_{G}\) is not injective.
Note that $\varphi$ is ``injective mod $G$'', meaning that $\varphi(z_1)=\varphi(z_2)$
if and only if $z_1=g\cdot z_2$ for some $g\in G.$
%\end{remark}

\bigskip
\noindent \textbf{The Symplectic Normal Space and the Witt-Artin Decomposition}
Let $G$ be a Lie group acting symplectically and properly on $\left(M,\omega\right),$
and let $z\in M.$ 
The Witt-Artin decomposition 
is a splitting
\begin{align}
\label{E:Witt}
T_zM=T_1\oplus\left(T_0\oplus N_0\right)\oplus N_1
\end{align}
such that 
$T_1\oplus T_0 = \mathfrak{g}\cdot z$
(so ``T'' is for ``tangent'' and ``N'' is for ``normal''),
and
$N_1\oplus T_0 = \left(\mathfrak{g}\cdot z\right)^\omega$
(the symplectic complement to $ \mathfrak{g}\cdot z$)
and 
each of the three spaces $T_1, \left(T_0\oplus N_0\right)$ and $N_1$
is a symplectic subspace of $\left(T_z M,\omega(z)\right).$
The decomposition can be chosen to be $G_z$-invariant,
where $G_z$ is the isotropy group of $z.$

We now define the decomposition, which is not unique, due to a choice of complements.
Though it can be defined more generally, we
will assume the existence of an $\Ad^*$-equivariant momentum map $J.$
Let $\mu=J(z)$
and let $G_{\mu}$ be the isotropy group
of $\mu$ under the coadjoint action. 
Let $\mathfrak{g}$ and $\mathfrak{g}_\mu $be the Lie algebras of $G$ and $G_\mu.$
The well-known {}``Reduction Lemma'' states that $\left(\mathfrak{g}\cdot z\right)^{\omega}=\ker dJ(z)$
and $\mathfrak{g}\cdot z\cap\left(\mathfrak{g}\cdot z\right)^{\omega}=\mathfrak{g}_{\mu}\cdot z.$
We define $T_0=\mathfrak{g}_{\mu}\cdot z.$
Let $T_1$ be a $G_z$-invariant complement to $T_0.$
which always exists since $G_z$ is compact.
Similarly, let $N_1$ be a $G_z$-invariant complement to $\mathfrak{g}_{\mu}\cdot z$ in $\ker dJ(z)$.
Since the kernel of $\omega$ restricted to 
either $\mathfrak{g}\cdot z$ or $\ker dJ(z)$ is $\mathfrak{g}_{\mu}\cdot z,$
the subspaces $T_1$ and $N_1$ are symplectic.
It can be shown that $T_0$ is a Lagrangian subspace of 
$\left(T_1\oplus N_1\right)^\omega$ and
there exists a Lagrangian subspace $N_0$ of $\left(T_1\oplus N_1\right)^\omega$ such that
$\left(T_1\oplus N_1\right)^\omega = T_0 \oplus N_0$
is a $G_z$-invariant splitting.
Then it can be shown that the Witt-Artin decomposition (Equation \ref{E:Witt}) has
the properties stated above, and that there is a $G_z$-invariant isomorphism of $N_0$
with $\left(\mathfrak{g_\mu}/\mathfrak{g}_z\right)^*.$
See \cite{OR04} for proofs and further discussion.

The space $N_1$ is often called the \emph{symplectic normal space}
to the group orbit through $z$,
but we will reserve this term for the isomorphic space
\[
N_s=\ker dJ(z)/\mathfrak{g}_{\mu}\cdot z,
\]
which has a reduced symplectic bilinear form defined as in Theorem \ref{linred}.
Since $G_z\subset G_\mu,$ it is easy to show that the tangent-lifted action 
of $G_z$ on $T_{z}M$ leaves $\ker dJ(z)$ invariant
and descends to a symplectic action on $N_{s}$ given by 
\begin{align} \label{E:HonNs}
h\cdot \left(v+\mathfrak{g}_{\mu}\cdot z\right) = T\Phi_h(v)+\mathfrak{g}_{\mu}\cdot z,
\end{align}
and that $N_s$ is $G_z$-equivariantly symplectomorphic to $N_1$ for any choice of $N_1$ as above.
Note that in the case of free actions, $N_s$ is symplectomorphic to the tangent space at 
$\left[ z\right]_{G_{\mu}}$ to the 
reduced space $J^{-1}\left(\mu\right)/G_{\mu}.$

%The tangent-lifted action of $H$ on $T_{z}M$ is a symplectic and
%leaves $\mathfrak{g}\cdot z$ invariant, so it also leaves $\left(\mathfrak{g}\cdot z\right)^{\omega}$
%invariant; it follows that the $H$ action restricts to an action
%on $\left(\mathfrak{g}\cdot z\right)^{\omega},$ which in turn and
%descends to one on $N_{s}.$ The action is clearly symplectic. We
%have shown: 

\bigskip
\noindent \textbf{Slice Theorems}
While symplectic reduction can be generalised to singular momentum values,
the resulting spaces are in general not smooth \cite{SL91,BL97,OR04}.
It thus makes sense to consider a related question: how far can we factor out a symmetry without 
losing smoothness? Slice theorems are local answers to this question
(or ``semilocal'', since the model spaces are ``global in the group direction'').
They say that such a space is locally isomorphic
to a twisted product of the group and a "slice" transverse to the group orbit.
In the case of a free action, there is no twisting; in general, the twisting is by the 
isotropy group of the point at which the local model is based.

We first define twisted products, slices and tubes and state a version of Palais' slice theorem for
the category of Lie group symmetries on manifolds. We then 
define the symplectic normal space and state 
the Hamiltonian Slice Theorem
of Marle, Guillemin and Sternberg.
Our new cotangent bundle slice theorem will be presented later as Theorem \ref{mainGmuG}.

Let $H$ be a Lie subgroup of a Lie group $G,$ and $N$ is a manifold
on which $H$ acts. Consider
the following two left
actions on $G$ $\times N:$\begin{align}
K\text{ acts by \emph{twisting}} & :k\cdot\left(g,n\right)=\left(gk^{-1},k\cdot n\right)\label{E:defGKact}\\
G\text{ acts by \emph{left multiplication}} & :g^{\prime}\cdot\left(g,n\right)=\left(g^{\prime}g,n\right).\nonumber \end{align}
 These actions are easily seen to be free and proper.
The \emph{twisted product} $G\times_{H}N$
is the quotient of $G\times N$ by the \emph{twist action} 
$h\cdot\left(g,n\right)=\left(gh^{-1},h\cdot n\right).$
It is a smooth manifold (since the twist action is free); in fact
$G\times_{H}N\rightarrow G/H$ is the vector bundle associated to
the $H$ action on $N.$ The \emph{left multiplication action} of
$G$ on itself commutes with the twist action and drops to a smooth
$G$ action on $G\times_{H}N,$ namely 
$g^{\prime}\cdot\left[g,n\right]_{H}=\left[g^{\prime}g,n\right]_{H}.$

Now consider a $G$ action on a manifold $M,$ and a point $z\in M,$
and let $H=G_{z}$ be the isotropy subgroup of $z.$ A \emph{tube}
for the $G$ action at $z$ is a $G$-equivariant diffeomorphism from
some twisted product $G\times_{H}N$\ to an open neighbourhood of
$z$ in $M,$ that maps $\left[e,0\right]_{H}$ to $z.$ The space
$N$ may be embedded in $G\times_{H}N$ as $\left\{ \left[e,n\right]_{H}:n\in N\right\} ;$
the image of the latter by the tube is called a \emph{slice.} A \emph{slice
theorem} is a theorem that guarantees the existence of a tube under
certain conditions. Palais \cite{P61} was the first to prove a slice
theorem for proper actions. Many smooth versions of his original theorem
are in common use. A proof of the
following version appeared in an appendix to an earlier version of the present article, but this appendix
has now been moved to the author's website due to space considerations and 
because the proof has recently appeared (with permission) in \cite{OR04}.

\begin{theorem}
[``Palais' Slice Theorem'']\label{slice}Let \(G\) be a Lie group acting properly
and smoothly on a manifold \(M,\) and let \(z\in M.\) Let \(H=G_{z}\) be the
isotropy group of \(z,\) and let \(N\) be any \(H\)-invariant complement to
\(\mathfrak{g}\cdot z.\) Choose a local \(H\)-invariant Riemannian metric around \(z\)
(such a metric always exists), and let \(\exp_{z}\) be the corresponding
Riemannian exponential based at \(z.\) Then there exists an \(H\)-invariant
neighbourhood \(U\) of \(0\) in \(N\) such that the map
\begin{align*}
\tau:G\times_{H}U  &  \rightarrow M\\
\left[  g,n\right]  _{H}  &  \longmapsto g\cdot\exp_{z}n
\end{align*}
is a tube for the \(G\) action at \(z.\)
\end{theorem}

Because of this theorem, an \(H\)-invariant complement to \(\mathfrak{g}\cdot z\) is 
sometimes called a \emph{linear slice} to the $G$ action at $z.$
 
In the case of a linear action, we can replace ``\(\exp_{z}n\)'' with ``\(z+n\)''
in the above statement. Indeed, it is easy to prove that a constant
\(H\)-invariant Riemannian metric always exists; the corresponding exponential
will be the map \(n\mapsto z+n.\) Alternatively, a special version of the slice
theorem for linear actions can be proven directly; in doing so, one notes that
the map $\tau$ will be a tube whenever the neighbourhood $U$ is small enough that $\tau$
is injective. Thus we have the following result.

\begin{theorem} [Slice theorem for linear actions]
\label{linslice}
Let \(G\) be a Lie group acting properly, smoothly and linearly
on \(\mathbf{R}^{n}\) and let \(z\in\mathbf{R}^{n}.\) Let \(H=G_{z}\) be the
isotropy group of \(z,\) and let \(N\) be any \(H\)-invariant complement to
\(\mathfrak{g}\cdot z.\) 
Then there exists an \(H\)-invariant neighbourhood \(U\) of
\(0\) in \(N\) such that the map
\begin{align*}
\tau:G\times_{H}U  &  \rightarrow\mathbf{R}^{n}\\
\left[  g,n\right]  _{H}  &  \longmapsto g\cdot\left(  z+n\right)
\end{align*}
is injective. Given any such \(U,\) the map \(\tau\) is a tube for the \(G\) action
at \(z.\)
\end{theorem}

\bigskip
\noindent \textbf{The Hamiltonian Slice Theorem}
Now suppose that $G$ acts \emph{symplectically} on a symplectic manifold
$\left(M,\omega\right).$ We would like tube $\tau$ to be symplectic.
One could obviously pull back the symplectic form on $M$ by the
diffeomorphism given by Palais' slice theorem, 
but the resulting symplectic structure need not be
simple or ``natural''.
The Hamiltonian
slice theorem defines a ``natural'' symplectic form on a space $G\times_{H}N$
and then shows that this space is $G$-equivariantly locally symplectomorphic with $M.$

The Hamiltonian slice theorem was first proven by Marle \cite{Marl85} and Guillemin and Sternberg \cite{GS84},
for compact groups $G$
and extended to proper actions of arbitrary groups by Bates and Lerman \cite{BL97}.
The Hamiltonian
slice theorem is called ``Hamiltonian'' because it assumes that
the $G$ action is globally Hamiltonian, meaning that it has a
globally defined
$\Ad^{\ast}$\emph{-}equivariant
momentum map $J$.
We note that this assumption has
been removed by Ortega and Ratiu \cite{OR02b}
and Scheerer and Wulff \cite{SW01}.
The Hamiltonian
slice theorem is sufficiently general for the present article, since
all cotangent-lifted actions have an $\Ad^{\ast}$\emph{-}equivariant
momentum map. We now present the Hamiltonian
slice theorem, following closely the presentations
in \cite{Ort98} and \cite{SL91}, to which we refer the reader for details and proofs. 

We are assuming that $G$ acts symplectically and properly on $\left(M,\omega\right),$ 
with $\Ad^{\ast}$\emph{-}equivariant
momentum map $J$.
Let $z\in M$ and let $H=G_{z}$ be the isotropy group of $z.$
Let $\mu=J\left(z\right)$ and let $G_{\mu}$ be the isotropy group
of $\mu$ under the coadjoint action. Note that $H\subset G_{\mu},$
by the $\Ad^{\ast}$-equivariance of $J.$ Let $\mathfrak{h}$,
$\mathfrak{g}$ and $\mathfrak{g}_{\mu}$ be the Lie groups of $H,$
$G$ and $G_{\mu}$ respectively. 
Let $\mathfrak{m}$ be an $H$-invariant complement
to $\mathfrak{h}$ in $\mathfrak{g}_{\mu}.$
Recall from above that the \emph{symplectic normal space} at $z$ is 
$N_{s}=\ker T_{z}J\,/\,\mathfrak{g}_{\mu}\cdot z$
and that there is a reduced symplectic form and a natural $H$ action on $N_s.$
This $H$ action and the coadjoint action of $H$ on $\mathfrak{m}^{\ast}$
define an $H$ action on $\mathfrak{m}^{\ast}\oplus N_{s},$ allowing
us to define the twisted product 
\[
G\times_{H}\left(\mathfrak{m}^{\ast}\oplus N_{s}\right).
\]
This will be the model space of the Hamiltonian slice theorem.
Recall from above that there is an isomorphism
$\mathfrak{m}^{\ast}\cong N_0$ and a symplectomorphism
$N_{s}\cong N_1,$ both $H$-equivariant, where $N_0$ and $N_1$ are
components in the Witt decomposition. The sum $N_0\oplus N_1$ is  linear slice at $z.$
Thus the model space $G\times_{H}\left(\mathfrak{m}^{\ast}\oplus N_{s}\right)$
can be considered to be a special case of the model space $G\times_{H}N$
in Theorem \ref{slice}. 

We now define the symplectic form on the space $G\times_{H}\left(\mathfrak{m}^{\ast}\oplus N_{s}\right),$%
beginning with a presymplectic form (i.e., a closed two-form) on 
$G\times\mathfrak{g}_{\mu}^{\ast}\oplus N_{s}.$
First, let $\Omega_{c}$ be the pull-back to $G\times\mathfrak{g}_{\mu}^{\ast}$
of the canonical form on $T^{\ast}G$ by the map
% \begin{align*}
$\, G\times\mathfrak{g}  _{\mu}^{\ast}\to T^{\ast}G,\ 
\left(g,\nu\right) \mapsto TL_{g^{-1}}^{\ast}\nu.
$
%\end{align*}
 Second, let $\Omega_{\mu}$ be the pull-back by the map 
 %\begin{align*}
$\, G\times\mathfrak{g}_{\mu}^{\ast}  \to\mathcal{O}_{\mu},\
\left(g,\nu\right)  \mapsto Ad_{g^{-1}}^{\ast}\mu,
$
%\end{align*}
 of the KKS symplectic form $\omega_{\mathcal{O}_{\mu}}^{+}$
 (defined in Eq. \ref{E:KKS}).
 Third, let $\omega_{N_{s}}$ be the reduced symplectic bilinear form
on $N_{s}.$ The sum $\Omega_{Z}$ $=\Omega_{c}+\Omega_{\mu}+\Omega_{N_{s}}$
is a presymplectic form on $Z=G\times\left(\mathfrak{g}_{\mu}^{\ast}\oplus N_{s}\right)$. 

Consider the twist action \emph{}of $H$ on $Z$ corresponding to
the coadjoint action of $H$ on $\mathfrak{g}_{\mu}^{\ast}$ and the $H$ action on $N_{s}$
inherited from the lifted action of $G$ on $T_zM.$
 The $H$-action on $N_{s}$ has an
$H$-equivariant momentum map $J_{N_{s}}$
(as does any linear symplectic action). 
%given by $\left\langle J_{N_{s}}\left(v\right),\eta\right\rangle =\frac{1}{2}\omega_{N_{s}}\left(\eta\cdot v,v\right)$
%(this is a standard formula for linear actions). 
One can check
that the twist action of $H$ on $Z$ is globally Hamiltonian with
respect to $\Omega_{Z},$ with momentum map $J_{H}:\left(g,\sigma,v\right)\longmapsto J_{N_{s}}\left(v\right)-\left.\sigma\right|_{\mathfrak{h}}.$
If we identify $\mathfrak{m}^*$ with $\mathfrak{k}^\circ\subset \mathfrak{g}_\mu^*,$ then
the following map is well-defined,
\begin{align}
l:G\times\mathfrak{m}^{\ast}\oplus N_{s} & \longrightarrow J_{H}^{-1}\left(0\right)\subset Z\label{E:l}\\
\left(g,\sigma,v\right) & \longmapsto\left(g,\sigma+J_{N_{s}}\left(v\right),v\right)\nonumber \end{align}
It is clearly an $H$-equivariant
diffeomorphism (with respect to twist action of $H$ on $G\times\mathfrak{m}^{\ast}\oplus N_{s}$
defined earlier). This map descends to a diffeomorphism $L$ defined by the following
commutative diagram, \begin{equation}
\begin{array}{ccc}
G\times\left(\mathfrak{m}^{\ast}\oplus N_{s}\right) & \overset{l}{\longrightarrow} & J_{H}^{-1}\left(0\right)\subset Z=G\times\left(\mathfrak{g}_{\mu}^{\ast}\oplus N_{s}\right)\\
\downarrow\pi_{H} &  & \downarrow\pi_{Z,H}\\
Y=G\times_{H}\left(\mathfrak{m}^{\ast}\oplus N_{s}\right) & \underset{\cong}{\overset{L}{\longrightarrow}} & J_{H}^{-1}\left(0\right)/H,\end{array}\label{E:L}\end{equation}
 where $\pi_{H}$ and $\pi_{Z,H}$ are the obvious projections. 

We define the presymplectic form $\omega_{Y}$ on $G\times_{H}\left(\mathfrak{m}^{\ast}\oplus N_{s}\right)$
as the pull-back by $L$ of the reduced presymplectic form on $J_{H}^{-1}\left(0\right)/H$
corresponding to $\Omega_{Z}.$ It can be shown that there exists
a $G$-invariant neighbourhood of $\left[e,0,0\right]_{H}$ in $Y$
in which $\omega_{Y}$ is symplectic. 

Finally, note that there is
left $G$-action on $Y$ given by $
g^{\prime}\cdot\left[g,\sigma,v\right]_{H}=\left[g^{\prime}g,\sigma,v\right]_{H}$
It is easy to check that this is symplectic %
with respect to $\omega_{Y}.$
We can now state the Hamiltonian Slice Theorem, also known as the
Marle-Guillemin-Sternberg Normal Form.
 
\begin{theorem}
[Hamiltonian Slice Theorem]\label{marle}In
the above context, there exists a symplectic tube from
%\(G\)-invariant neighbourhood of \(\left[  e,0,0\right]  _{H}\) in the manifold
\(Y=G\times_{H}\left(  \mathfrak{m}^{\ast}\oplus N_{s}\right)  \) to $M$
that maps 
%\(G\)-invariant neighbourhood of \(z\) in \(M,\) such that 
\(\left[  e,0,0\right]_{H}\) to \(z.\)
\end{theorem}

It can be shown that the momentum map 
of the $G$ action on $Y$
is 
\[
J_Y\left(\left[g,\sigma,v\right]_H\right) = \Ad^*_{g^{-1}}\left( \mu + \sigma +  J_{N_s}(v)\right)
\]

\section{Commuting reduction\label{sectcommred}}

In this section we consider a manifold with two commuting symplectic
actions. We first review regular and singular commuting reduction and then introduce
a new ``tangent-level version'' of commuting reduction, which we will use in the next section 
in our analysis of the symplectic normal space of a cotangent-lifted action.

We have already seen an example of commuting symplectic actions
in the presentation of the Hamiltonian slice
theorem: the $G$ and $H$ actions on the manifold $G\times\left(\mathfrak{g}_{\mu}^{\ast}\oplus N_{s}\right)$ (see Equation \ref{E:L}).
In this context, commuting reduction leads to 
a singular local normal form for a symplectic reduced space, Theorem \ref{SLNFred}.
A second example of commuting symplectic actions, key to the rest
of this article, will appear in the next section: a bundle $T^{\ast}\left(G\times A\right)$
with the cotangent lifts of the left
multiplication action of $G$ and the twist action of a subgroup $K$ of $G.$ Commuting
reduction in this context leads to a cotangent-bundle specific local normal form
for a symplectic reduced space, Theorem
\ref{mainGK}.``Tangent-level
reduction'' in this context will be used to characterise the symplectic
normal space of a cotangent-lifted action: see Theorem \ref{mainVm} and following results.

Let \(G\) and \(K\) be Lie groups acting symplectically and
properly on a symplectic manifold \(M,\) with equivariant momentum
maps \(J_{G}\) and \(J_{K}\) respectively, and suppose that the actions commute.
Let \(\mu \in\mathfrak{g}^{\ast}\) and \(\nu \in\mathfrak{k}^{\ast}\).
The idea of commuting reduction is to first reduce by the $K$ action
(say) and then reduce the $K$-reduced space by the induced $G$ action;
and then switch the order, reducing first by $G$ and then by $K.$
Under very general conditions, the two doubly-reduced spaces are isomorphic.
We first state the ``regular version'' of commuting reduction,
due to Marsden and Weinstein \cite{MW74}; the key assumption here
is that all of the group actions are free. 

\begin{theorem}
[Regular Commuting Reduction]\label{RCR}
In the above context, suppose that $G$ and $K$ act freely and \(J_{K}\) is \(G\)-invariant
and \(J_{G}\) is \(K\)-invariant. Then \(G\) induces a symplectic action on \(M_{\nu
}:=J_{K}^{-1}\left(  \nu\right)  /K_{\nu}\) with equivariant momentum map
\(J_{\bar{G}}\) determined by \(J_{\bar{G}}\circ\pi_{K_{\nu}}=J_{G}\) (where
\(\pi_{K_{\nu}}:\) \(M\rightarrow M/K_{\nu}\) is projection, and both sides of the
equation are restricted to \(J_{K}^{-1}\left(  \nu\right)  \)). If the reduced
\(G\) action is free, then the reduced space for this action at \(\mu\) is
symplectomorphic to the reduction of \(M\) at \(\left(  \mu,\nu\right)  \) by the
product action of \(G\times K.\)
%\marginpar{are all the conditions needed?
%\(J_{K}?\) also, free isn't needed, just constant orbit type?}
\end{theorem}

Note that applying this theorem a second time, with the roles of $G$ and $K$ reversed, 
shows that the reduced
space at \(\nu\) for the action of \(K\) on \(J_{G}^{-1}\left(  \mu\right)
/G_{\mu}\) is symplectomorphic to the reduced space at \(\mu\) for the
action of \(G\) on \(M_{\nu}\).

Sjamaar and Lerman \cite{SL91}, working with reduction at zero of
compact group actions, showed that 
a similar result holds even if the actions
are not free. In this case, the reduced spaces need not
be smooth manifolds, but are Poisson varieties. In the general case,
for proper actions and arbitrary momentum values, we need to add the
hypotheses that $G_{\mu}$ and $K_{\nu}$ are compact and that $\mathcal{O}_{\mu}$
and $\mathcal{O}_{\nu}$ are locally closed, the latter for reasons
discussed in \cite{MMOPR}. 

\begin{theorem}
[Singular Commuting Reduction]\label{SCR}
In the above context, suppose that \(J_{K}\) is \(G\)-invariant, \(J_{G}\) is \(K\)-invariant,
 \(G_{\mu}\) and \(K_{\nu}\) are compact and the coadjoint
orbits \(\mathcal{O}_{\mu}\) and \(\mathcal{O}_{\nu}\) are locally closed. Then
\(G\) induces a Poisson action on \(M_{\nu}=J_{K}^{-1}\left(  \nu\right)
/K_{\nu},\) with equivariant momentum map \(J_{\bar{G}}\) determined by
\(J_{\bar{G}}\circ\pi_{K_{\nu}}=J_{G}.\) The reduced space for the action of \(G\)
on \(M_{\nu}\) at \(\mu\) is Poisson diffeomorphic to the reduction of \(M\) at
\(\left(  \mu,\nu\right)  \) by the product action of \(G\times K.\)
\end{theorem}

It follows that the reduced
space at \(\nu\) for the action of \(K\) on \(J_{G}^{-1}\left(  \mu\right)
/G_{\mu}\) is Poisson diffeomorphic to the reduced space at \(\mu\) for the
action of \(G\) on \(M_{\nu}\).

The Hamiltonian Slice Theorem (Theorem \ref{marle}), together with 
singular commuting reduction, applied to the $G$ and $H$ actions on  $G\times\left(\mathfrak{g}_{\mu}^{\ast}\oplus N_{s}\right)$ (see Equation \ref{E:L}),
can be used to deduce the following local normal form for a symplectic reduced space
(when $G_{\mu}$ is compact). 
The result was first published by Sjamaar and Lerman
\cite{SL91} for $\mu=0;$ the general case is due to Bates and Lerman
\cite{BL97}. The proof  given in \cite{BL97} does not use a commuting reduction theorem
and does not require $G_{\mu}$ compact.

\begin{theorem} \label{SLNFred}
Let \(G\) act properly on the symplectic manifold \(\left(  M,\omega\right)  \)
with equivariant momentum map \(J.\) Let \(z\in M\) and \(H=G_{z}\) and
\(\mu=J\left(  z\right)  ,\) and let \(N_{s}\) be the symplectic normal space to
\(\mathfrak{g}\cdot z.\) Assume that the coadjoint orbit \(\mathcal{O}_{\mu}\) is
locally closed. Then there is a local Poisson diffeomorphism between the
reduced space \(J^{-1}\left(  \mu\right)  /G_{\mu}\) and the reduced space at
\(0\) for the \(H\) action on \(N_{s}.\)
\end{theorem}

In the case of cotangent-lifted actions, our analysis of the symplectic normal space,
in the next section,
together with the above theorem comprise 
a cotangent-bundle-specific local normal form for symplectic reduced spaces, as
we note later  in Remark \ref{CBSLNFred}.

We now introduce another approach to singular commuting reduction,
assuming that the original actions are free but not assuming that the quotient action on the
once-reduced space is free.
Recall that, in the case of a free action, the symplectic normal space 
``is'' the tangent space to the reduced space. 
This observation suggests studying symplectic normal spaces in place of the
possibly singular  doubly-reduced spaces.

%Though this
%does not hold in the general for non-free actions
%
%\footnote{In the case of non-free actions, the reduced spaces may of course
%be singular. Singular symplectic reduction may be applied (see \ref{BL97}),
%giving a stratification into smooth symplectic pieces, but the tangent
%spaces to these pieces are not necessarily isomorphic to $N_{s}\left(z\right)$
%(consider for example the cotangent lift of the standard $SO\left(2\right)$
%action on $\mathbf{R}^{2},$ with $z=\left(\mathbf{0},\mathbf{0}\right)).$%
%},

Since symplectic normal spaces are quotients, the following lemma and notation will be useful;
the lemma is easily checked.

\begin{lemma} \label{fbar}
Let $\omega_A$ and $\omega_B$
be bilinear forms on vector spaces $A$ and $B$, respectively.
Suppose \(f:A\rightarrow B\) satisfies $f^*\omega_B = \omega_A.$
Then the quotient map
\(\bar{f}:A/\ker\left(  \omega_A\right)  \rightarrow B/\ker\left(  \omega_B\right) \) 
is well-defined and
injective. If \(f\) is surjective, then \(\bar{f}\) is bijective. 
If $\omega_A$ and $\omega_B$ are presymplectic (i.e. skew-symmetric)
then $\bar{f}$ is symplectic.
Also, if \(\bar{g}:B/\ker\left(  \sigma_B\right)  \rightarrow C/\ker\left(  \sigma_C\right) \) 
is defined similarly then \(\overline{f\circ g}=\overline
{f}\circ\overline{g}.\)
\end{lemma}

\begin{theorem}
[ ``Tangent-level'' commuting reduction]\label{TLCR}Let \(G\) and \(K\) be free,
symplectic, commuting actions on a symplectic manifold \(M,\) with momentum maps
\(J_{G}\) and \(J_{K}\) respectively. Then the product action of \(G\times K\) has
momentum map given by \(J_{G\times K}\left(  x\right)  =\left(  J_{G}\left(
x\right)  ,J_{K}\left(  x\right)  \right)  .\) Let \(x\in M\)  and \(\left(
\mu,\nu\right)  =J_{G\times K}\left(  x\right)  .\) The symplectic normal space
at \(x\) for the product action is\[
N_{s}\left(  x\right)  =\ker T_{x}J_{G\times K}/\left(  \mathfrak{g}_{\mu}\cdot
x+\mathfrak{k}_{\nu}\cdot x\right)  .
\]
Suppose further that \(G\) acts properly and that \(J_{G}\) is \(Ad^{\ast}\)-equivariant and that \(J_{G}^{-1}\left(  \mu\right)  \) is \(K\)-invariant. Then
the quotient action of \(K\) on \(J_{G}^{-1}\left(  \mu\right)  /G_{\mu}\) is
symplectic with respect to the reduced symplectic form, and its momentum map
\(J_{\bar{K}}\) satisfies \(J_{\bar{K}}\circ\pi_{G_{\mu}}=\left.
J_{K}\right|  _{J_{G}^{-1}\left(  \mu\right)  }\) (where \(\pi_{G_{\mu}}:J_{G}^{-1}\left(  \mu\right)  \rightarrow J_{G}^{-1}\left(  \mu\right)
/G_{\mu}\) is projection). The map \(\left(  g,k\right)  \longmapsto k\) is a Lie
group isomorphism from \(\left(  G\times K\right)  _{x}\) to \(K_{\left[
x\right]  _{G_{\mu}}}\) (where \(\left[  x\right]  _{G_{\mu}}=G_{\mu}\cdot x).\)
We identify these two groups and call them \(H.\) Let \(N_{s}\left(  \left[
x\right]  _{G_{\mu}}\right)  \) be the symplectic normal space at \(\left[
x\right]  _{G_{\mu}}\) for the \(K\) action on \(J_{G}^{-1}\left(  \mu\right)
/G_{\mu}.\) Let \(H\) act on each symplectic normal space, as in Equation \ref{E:HonNs}.
Then the following is an \(H\)-equivariant vector space
symplectomorphism,
\begin{align*}
\overline{T_{x}\pi_{G_{\mu}}}:N_{s}\left(  x\right)   &  \rightarrow
N_{s}\left(  \left[  x\right]  _{G_{\mu}}\right)  \\
v+\left(  \mathfrak{g}_{\mu}\cdot x+\mathfrak{k}_{\nu}\cdot x\right)   &
\longmapsto T\pi_{G_{\mu}}\left(  v\right)  +\left(  \mathfrak{k}_{\nu}\cdot\left[  x\right]  _{G_{\mu}}\right)  .
\end{align*}
\end{theorem}

\begin{proof}
It is easily verified that the product action has the given momentum map.
Since  \(\left(  G\times K\right)  _{\left(  \mu
,\nu\right)  }=G_{\mu}\times K_{\nu}\) and the actions commute, 
we have \(\left(  \mathfrak{g}\oplus\mathfrak{k}\right)  _{\left(  \mu,\nu\right)  }\cdot x=\left(
\mathfrak{g}_{\mu}\oplus\mathfrak{k}_{\nu}\right)  \cdot x=\left(  \mathfrak{g}_{\mu}\cdot x+\mathfrak{k}_{\nu}\cdot x\right)  ,\) so the symplectic normal
space at \(x\) is \(N_{s}\left(  x\right)  =\ker T_{x}J_{G\times K}/\left(
\mathfrak{g}_{\mu}\cdot x+\mathfrak{k}_{\nu}\cdot x\right)  .\) The claims about
the quotient action of \(K\) on \(J_{G}^{-1}\left(  \mu\right)  /G_{\mu}\) are
part of regular commuting reduction (Theorem \ref{RCR}), and in any case are
easy to prove by ``diagram-chasing''.

We will now show that \(\theta:\left(  G\times K\right)  _{x}\rightarrow
K_{\left[  x\right]  _{G_{\mu}}},\left(  g,k\right)  \longmapsto k,\) is an
isomorphism. To show it's well-defined, let \(\left(  g,k\right)  \in\left(
G\times K\right)  _{x},\)
so $k\cdot x = g^{-1}\cdot x.$
 Since \(J_{G}\) is \(Ad^{\ast}\)-equivariant and
\(J_{G}^{-1}\left(  \mu\right)  \) is \(K\)-invariant, we have \(\mu=J_{G}\left(
x\right)  =J_{G}\left(  \left(  g,k\right)  \cdot x\right)  =g\cdot\mu,\) so
\(g\in G_{\mu}.\) This implies that \(k\in K_{\left[  x\right]  _{G_{\mu}}}.\) So
\(\theta\) is well-defined. It is clearly smooth, and a homomorphism. For every
\(k\in K_{\left[  x\right]  _{G_{\mu}}},\) we have $k\cdot x\in G_\mu x;$
since $G$ acts freely, there is a unique element $\gamma(k)\in G_\mu$ 
such that \(k\cdot x=\gamma\left(  k\right)
^{-1}\cdot x.\)
Clearly \(\left(  \gamma\left(  k\right)
,k\right)  \cdot x=x,\) so the map \(k\mapsto\left(  \gamma\left(  k\right)
,k\right)  \) is an inverse for \(\theta.\) The smoothness of \(\theta^{-1}\) is a
consequence of the implicit function theorem applied to the restricted action
\(F:G_{\mu}\times K_{\left[  x\right]  _{G_{\mu}}}\rightarrow G_{\mu}\cdot x\)
given by \(F\left(  g,k\right)  =\left(  g,k\right)  \cdot x.\) Indeed, note
that \(\left(  G\times K\right)  _{x}=F^{-1}\left(  x\right)  ,\) and that
\(D_{1}F\left(  g,k\right)  \) is surjective for every \(\left(  g,k\right)  \in
G_{\mu}\times K_{\left[  x\right]  _{G_{\mu}}},\) since the \(G\) action is free.
Hence \(\theta\) is a a Lie group isomorphism. We identify \(\left(  G\times
K\right)  _{x}\) with \(K_{\left[  x\right]  _{G_{\mu}}}\) via \(\theta,\) calling
both groups \(H.\)

Next, observe that 
\begin{align*}
\ker T_{x}J_{G\times K} &  =\ker T_{x}J_{G}\cap\ker T_{x}J_{K}=T_{x}J_{G}^{-1}\left(  \mu\right)  \cap\ker T_{x}J_{K}=\ker T_{x}\left(  \left.
J_{K}\right|  _{J_{G}^{-1}\left(  \mu\right)  }\right)  \\
&  =\ker T_{x}\left(  J_{\overline{K}}\circ\pi_{G_{\mu}}\right)  =\left(
T_{x}\pi_{G_{\mu}}\right)  ^{-1}\left(  \ker T_{\left[  x\right]  _{G_{\mu}}}J_{\bar{K}}\right).
\end{align*}
Since \(T_{x}\pi_{G_{\mu}}\) is surjective, this implies that \(T_{x}\pi_{G_{\mu}}\left(  \ker T_{x}J_{G\times K}\right)  =\ker T_{\left[
x\right]  _{G_{\mu}}}J_{\bar{K}}.\) 
The map $T\pi_{G_\mu}$ is a presymplectic submersion, by definition of the reduced symplectic form
on $J_G^{-1}(\mu)/G_\mu.$ 
Hence Lemma \ref{fbar} implies that 
\(\overline{T_{x}\pi_{G_{\mu}}},\) as
defined in the statement of the theorem, is a well-defined symplectic isomorphism
from \(N_{s}\left(  x\right)  \) to \(N_{s}\left(  \left[  x\right]  _{G_{\mu}}\right)  .\) 

The projection $\pi_{G_\mu}$ is $K$-equivariant, by definition of the quotient action of $K.$
Since we have already shown that \(\left(  g,k\right)  \in\left(  G\times K\right)  _{x}\)
implies \(g\in G_{\mu}\),
the \(H\)-equivariance of \(\pi_{G_{\mu}}\)
is easily checked.
The \(H\)-equivariance of  \(T_{x}\pi_{G_{\mu}}\), and hence \(\overline{T_{x}\pi_{G_{\mu}}}\), follows.
\end{proof}

\section{The symplectic normal space of a cotangent-lifted action\label{sectNs}}

The main result of the section will be a characterisation of the symplectic normal space $N_s$ to the
orbit of a cotangent-lifted action, given in Theorem \ref{mainVm}.
In two special cases this leads to new splittings of $N_s,$
given in Corollaries \ref{NsKsubGmu} and \ref{Nsalpha0}.
Our analysis of the special case
$G_q\subset G_{\mu}$ and much of the general set-up developed in this section
will be used later in the cotangent
bundle slice theorem (Theorem \ref{mainGmuG}). 
We also note implications for singular reduction, in Theorem \ref{mainGK} and 
Remark \ref{CBSLNFred},
and the reconstruction equations (bundle equations), in Equations 
\ref{E:reconGmuG} and \ref{E:reconalpha0}.

Let $G$ act properly by cotangent lifts on $T^{\ast}Q,$ with momentum
map $J,$ and let $z\in T_{q}^{\ast}Q$ and $\mu=J\left(z\right).$
Let $K=G_{q}$ and $H=G_z,$
and let $\mathfrak{g}, \mathfrak{g}_{\mu}, \mathfrak{k}, \mathfrak{h}$ be the Lie algebras
of $G,G_{\mu},K$ and $H.$

\begin{lemma} \label{cbiso} (i) $H\subset K$
(ii) $H\subset G_{\mu}$
(iii) $\mathfrak{k} \subset \ker \mu$
(iv) If \(K\) is normal in \(G,\) then \(K\subset G_{\mu}.\)
\end{lemma}

\begin{proof} Claim (i) is clear from $z\in T_q^*Q;$ (ii) follows from 
the equivariance of $J.$ (iii) The definition of \(J\) gives \(\left\langle \mu,\xi\right\rangle \)
\(=\left\langle z,\xi_{Q}\left(  q\right)  \right\rangle =0\) for all \(\xi \in \mathfrak{k}.\)
(iv) For every \(g\in G\) and \(k\in K\) we have
\(gkg^{-1}k^{-1}\in K.\) Differentiating with respect to \(g\) gives
\(\xi-Ad_{k^{-1}}\xi\in\mathfrak{k.}\) Thus, for every \(k\in K\) and \(\xi
\in\mathfrak{g},\) we have \(\left\langle Ad_{k^{-1}}^{\ast}\mu-\mu,\xi
\right\rangle =\left\langle \mu,Ad_{k^{-1}}\xi-\xi\right\rangle =0,\) in other
words \(k\in G_{\mu}.\)
\end{proof}

There exist simple examples in which $H$ is a proper subset of $K\cap G_{\mu}.$
The complex relationship between the different isotropy subgroups
is one of the key difficulties of the subject. 

We begin by applying Palais' slice theorem (Theorem \ref{slice}) to 
the configuration space $Q.$
Choose a $K$-invariant Riemannian metric on some
neighbourhood of $q$ in $Q$, and let $A$ be the orthogonal complement
to $\mathfrak{g}\cdot q$ in $T_{q}Q,$ written $A=$ $\left(\mathfrak{g}\cdot q\right)^{\perp}.$
By Palais' slice theorem (Theorem \ref{slice}), there exists a $K$-invariant
neighbourhood $U$ of $0$ in $A$ such that the map 
\begin{align} \label{E:defs} 
s:G\times_{K}U & \longrightarrow Q\\
\left[g,a\right]_{K} & \longmapsto g\cdot\exp_{q}a
\end{align}
is a $G$-equivariant embedding. The cotangent lift 
\begin{align} \label{E:Qslice}
T^{\ast}s^{-1}:T^{\ast}\left(G\times_{K}U\right)\rightarrow T^{\ast}Q
\end{align}
 is a $G$-equivariant symplectic embedding onto a neighbourhood of
$z$ (symplectic with respect to the standard cotangent bundle symplectic
forms). 

We next make the following key observation (explained fully in Proposition \ref{unrollK}):
\begin{quote}
\emph{
$T^{\ast}\left(G\times_{K}U\right)$
is a reduced space
for the lifted twist action of $K$
on $T^{\ast}\left(G\times U\right).$}
\end{quote}
We will apply commuting reduction to $T^{\ast}\left(G\times U\right),$ with the second action
being the lift of left multiplication by $G.$
We first fill in the details of the passage
to $T^{\ast}\left(G\times U\right),$ and state some basic facts for
later use.

Let $N$ be any manifold on which $K$ acts
(we have in mind $N=U$ or $N=A,$ but the following facts are general).
Recall from the Equation \ref{E:defGKact} the following two left
actions on $G\times N,$
which commute and are both free and proper:
\begin{align}
K\text{ acts by \emph{twisting}} & :k^K\cdot\left(g,n\right)=\left(gk^{-1},k\cdot n\right)\label{E:defGKact2}\\
G\text{ acts by \emph{left multiplication}} & :h^G\cdot\left(g,n\right)=\left(hg,n\right).\nonumber \end{align}
Note that, since $K$ is a subset of $G,$
there is room for confusion of the two actions, so we have introduced
superscripts to identify them.
Each of these actions has a corresponding tangent-lifted action on
$T\left(  G\times N\right)\cong TG\times TN$
and cotangent-lifted action on
$T^{\ast}\left(  G\times N\right)\cong T^*G\times T^*N.$
It is easy to see that these actions commute and are free and proper.

Throughout this article, we will identify $TG$ with $G\times\mathfrak{g}$ and
$T^*G$ with $G\times\mathfrak{g}^*$ by
left trivialisation,
\[
\begin{tabular}[c]{ll}$TG\overset{\cong}{\longrightarrow}G\times\mathfrak{g}\quad$and
&$T^*G\overset{\cong}{\longrightarrow}G\times\mathfrak{g}^*$ \\
$  TL_{g}(\xi)\longmapsto \left(  g,\xi\right)$
&$T^*L_{g^{-1}}(\nu) \longmapsto \left(  g,\nu\right) $
\end{tabular}
\]
where $L_{g}$ is left multiplication by $g.$ 
The following basic properties of the left and right multiplication actions are well known.
%(see for example \cite{SL91}). 

\begin{lemma}
\label{LRbasic}Let \(G\) be a Lie group. With
respect to the left trivialisations of \(TG\) and \(T^{\ast}G,\) the left and right 
multiplication actions of \(G\) on itself have the following lifted actions and
infinitesimal lifted actions:
\begin{align*}
\text{\emph{tangent}}\emph{:}\text{ }  &  h^{L}\cdot\left(
g,\xi\right)  =\left(  hg,\xi\right) 
&& h^{R}\cdot\left(  g,\xi\right)  =\left(
gh^{-1},\Ad_{h}\xi\right)
\\
\text{\emph{cotangent}}\emph{:}\text{ }  &  h^{L}\cdot\left(
g,\nu\right)  =\left(  hg,\nu\right) 
&&h^{R}\cdot\left(
g,\nu\right)  =\left(  gh^{-1},\Ad\nolimits_{h^{-1}}^{\ast}\nu\right)
\\
\text{\emph{infinitesimal tangent: }}  &  \eta^{L}\cdot\left(
g,\xi\right)  =\left(  \Ad_{g^{-1}}\eta,0\right)  \text{\emph{ }}
&&\eta^{R}\cdot\left(
g,\xi\right)  =\text{\emph{ }}\left(  -\eta,\Ad_{\eta}\xi\right) 
\\
\text{\emph{infinitesimal cotangent: }}  &  \eta^{L}\cdot\left(
g,\nu\right)  =\left(  \Ad_{g^{-1}}\eta,0\right)  \text{\emph{ }}
&&\eta^{R}\cdot\left(  g,\nu\right)  =\left(  -\eta
,-\Ad_{\eta}^{\ast}\nu\right)
\end{align*}
The cotangent-lifted actions have the
following momentum maps,
with respect to the canonical symplectic form on $T^*G:$
\begin{align*}
 &  J_{L}\left(
g,\nu\right)  =\Ad{}_{g^{-1}}^{\ast}\nu,
&&J_{R}\left(
g,\nu\right)  =-\nu.
\end{align*}
The momentum map  \(J_{L}\) is invariant under the right multiplication action,
and \(J_{R}\) is invariant under the left multiplication action.
\end{lemma}

There are obvious corresponding properties for the $G$ and $K$ actions on $G\times N.$
In particular, we have the following:

\begin{remark}\label{GKbasic}
Let $G$ and $K$ act on $G\times N$ as in Equation \ref{E:defGKact2}.
Then the momentum maps for the cotangent-lifted actions on 
$T^*(G\times N) \cong G\times \mathfrak{g}^* \times T^*N,$
with respect to the canonical symplectic form on $T^*(G\times N),$ are
\begin{align*}
 &  J_{G}\left(g,\nu,w\right)  =\Ad{}_{g^{-1}}^{\ast}\nu,
&&J_{K}\left(g,\nu,w\right)  =-\left.\nu\right|_\mathfrak{k}+J_N(w),
\end{align*}
where $J_N(w)$ is the momentum map for the cotangent-lifted action of $K$ on $T^*N.$
The previous lemma implies that
\(J_{G}\) is invariant under the twist action of \(K\) and
 \(J_{K}\) is invariant under the left multiplication action of \(G.\) 
If $N$ is a vector space, we can identify $TN$ with $N\times N$
and $T^*N$ with $N\times N^*,$ so
\[
\begin{tabular}[c]{ll}$T(G\times N)\cong
G\times\mathfrak{g}\times N\times N\quad$and
&$T^*\left(G\times N\right)\cong G\times\mathfrak{g}^*\times N\times N^*,$\end{tabular}
\]
where the first and third components
are the base space, and the second and fourth are the (co-)tangent fibers.
These identifications will be used throughout this article.
In these coordinates, and using the diamond notation (see Equation \ref{E:diamond}),
$
J_{K}\left(  g,\nu,a,\delta\right) =-\nu|_{\mathfrak{k}}+a\diamond_{\mathfrak{k}}\delta.
$
\end{remark}

We are now in a position to apply reduction theorems to the two actions on $T^*(G\times N).$
We begin by studying the reduced space at $0$ for the $K$ action, using
cotangent bundle reduction
(Theorem \ref{RPCBRZ}). Note that $(G\times N)/K = G\times_K N.$
The map $\varphi$ in Theorem \ref{RPCBRZ} takes the following form:
\begin{align} \label{E:cbphi2}
\varphi:\left(  J_{K}^{-1}\left(  0\right)  \subset T^{\ast}\left(  G\times
N\right)  \right)  &\rightarrow T^{\ast}\left(  G\times_{K}N\right), \quad
\left\langle \varphi\left(  p\right)  ,T\pi_{K}\left(  v\right)  \right\rangle
=\left\langle p,v\right\rangle,
\end{align}
where \(\pi_{K}:G\times N\rightarrow
G\times_{K}N\) is projection. 
Recall that $G$ has a quotient action on $G\times_K N,$ and so 
$G$ acts on $T^*(G\times_K N)$ by cotangent lifts.
The projection \(\pi_{K}\) is
\(G\)-equivariant by definition of the $G$ action on \(G\times_{K}A,\) so
\(T\pi_{K}\) is \(G\)-equivariant with respect to the tangent lifted actions,
from which it follows that $\varphi$ is $G$-equivariant.
Since  \(J_{K}^{-1}\left(  0\right)  \) is \(G\)-invariant, the $G$ action descends to one
on on \(J_{K}^{-1}\left(  0\right)  /K\). 
It is easily verified that this quotient action is symplectic;
in fact this claim is part of Theorem \ref{RCR} (regular commuting reduction).
Applying Theorems \ref{RPCBRZ} and \ref{RCR} gives the following result.

\begin{proposition}
\label{unrollK}
Let $G$ and $K$ act on \(T^{\ast}\left(  G\times N\right)  \) as above, with momentum maps
\(J_{K}\)  and $J_G$ respectively.
Let $\varphi$ be defined as in Equation \ref{E:cbphi2}.
Then \(\varphi\) is a \(G\)-equivariant \(K\)-invariant surjective submersion that
descends to a \(G\)-equivariant symplectomorphism\[
\bar{\varphi}:J_{K}^{-1}\left(  0\right)  /K\rightarrow T^{\ast}\left(
G\times_{K}N\right)  ,
\]
with respect to the reduced symplectic form on \(J_{K}^{-1}\left(  0\right)
/K\) and the canonical symplectic form on \(T^{\ast}\left(  G\times_{K}N\right)
.\) If  \(J'\) is the momentum map for the \(G\) action on
\(T^{\ast}\left(  G\times_{K}N\right)  ,\) then the restriction of \(J_{G}\) to
\(J_{K}^{-1}\left(  0\right)  \) equals \(J'\circ\varphi.\)
\end{proposition}

Taking $N=A$ and combining this result with the map $T^*s^{-1}$ from Equation \ref{E:Qslice},
we have a $G$-equivariant symplectic embedding,
\begin{equation}
\begin{array}{ccccc}
\left(\left(J_{K}^{-1}\left(0\right)\cap\left(G\times\mathfrak{g}^{\ast}\times U\times A^{\ast}\right)\right) 
\right) / K
 \overset{\bar{\varphi}}{\longrightarrow} & T^{\ast}\left(G\times_{K}U\right) & \overset{T^{\ast}s^{-1}}{\hookrightarrow} & T^{\ast}Q\\
\end{array}\label{E:phis}\end{equation}
In particular, there exists an $x\in J_{K}^{-1}\left(0\right)$ such that 
$T^*s^{-1}\left(\varphi(x)\right)
=z\in T_q^*Q.$
Since $s\left(\left[e,0\right]_{K}\right) = q,$  we see that $\varphi(x)$
has base point $\left[e,0\right]_{K}.$ Since $\varphi$
covers $\pi_{K}:G\times U\rightarrow G\times_{K}U$ and is $K$-invariant,
we can choose $x$ to have base point $\left(e,0\right);$ 
in fact, since $\varphi$ is injective, this uniquely determines $x.$
So $x=(e,\nu,0,\alpha),$
for some $\nu\in\mathfrak{g}^*$ and some $\alpha\in A^*.$
Using Proposition \ref{unrollK}, we have $\nu= J_G(x)= J'\left(\varphi(x)\right)=J(z) =\mu.$
We can also show that $\alpha = \left.z\right|_A.$ 
Indeed, for every \(v\in A,\) we have 
\begin{align*}
\left\langle \alpha,v\right\rangle  
& = \left\langle \varphi\left(  e,\mu,0,\alpha\right), T\pi_K(e,0,0,v)\rangle\right)
=\left\langle T^{\ast}s\,\left(
z\right)  ,T\pi_{K}\left(  e,0,0,v\right)  \right\rangle \\
&  =\left\langle z,T\left(  s\circ\pi_{K}\right)  \left(  e,0,0,v\right)
\right\rangle
=\left\langle z,v\right\rangle ,
\end{align*}
where in the last line we have used \(s\circ\pi_{K}\left(  g,a\right)
=g\cdot\exp_{z}\,a,\) and the fact that the derivative at \(0\) of \(\exp_{z}\) is
the identity. In summary, we have shown:

\begin{lemma}\label{xalpha}
Let $\alpha = \left.z\right|_A$ and $x=\left(e,\mu,0,\alpha\right).$ Then
$T^*s^{-1}\left(\varphi(x)\right)
=z.$
\end{lemma}

\begin{remark}\label{H}
Recall that $H=G_z$ and note that $G_z=G_{\left[x\right]_K},$ where $\left[x\right]_K= \pi_K(x).$ 
By definition of the $G$ and $K$ actions
on $T^*(G\times A),$ we have \(H=G_{\mu}\cap K_{\alpha}\).
We note for later use that if \(H=K\) then  \(K_{\alpha}=K;\)
while if $G_\mu=G$ then $H=K_\alpha.$
\end{remark}

Applying singular commuting reduction (Theorem \ref{SCR}) gives the following picture:
\[
\begin{array}{ccc}
 & T^{\ast}\left(G\times A\right)\\
\text{reduction by }G\swarrow &  & \searrow\text{reduction by }K\\
\\?\cong J_{G}^{-1}\left(\mu\right)/G_{\mu} &  & J_{K}^{-1}\left(0\right)/K\overset{\text{\emph{local}}}{\cong}T^{\ast}Q\\
\\\text{reduction by }K\searrow &  & \swarrow\text{reduction by }G\\
 & ?\overset{\text{\emph{local}}}{\cong}J^{-1}\left(\mu\right)/G_{\mu}\end{array}\]

We now compute $J_{G}^{-1}\left(\mu\right)/G_{\mu}.$ 
Note that the $G$ action leaves $A$ untouched, so 
\[
J_{G}^{-1}\left(\mu\right)/G_{\mu}\cong \left( J_{L}^{-1}\left(\mu\right)/G_{\mu}\right) \times T^*A,
\]
where $J_{L}^{-1}\left(\mu\right)/G_{\mu}$ is the symplectic reduced space at $\mu$
for the lifted left multiplication action of $G$ on $T^*G.$
It is well-known that
$ J_{L}^{-1}\left(\mu\right)/G_{\mu}$ is symplectomorphic to  $\mathcal{O}_{\mu}$
with the KKS symplectic form 
\(\omega_{\mathcal{O}_{\mu}}^{-}\)
(defined in Eq. \ref{E:KKS});
see for example, Appendix B.4 of \cite{CB97}.
The isomorphism is $\left[g,\nu\right]_{G_{\mu}}\mapsto \nu$
(using the left trivialisation of $T^*G$.)
We have almost proven the following proposition; the remain claims in it are easily verified.

\begin{proposition}
\label{psi}The map \(\theta\) defined by\begin{align}
\theta:J_{G}^{-1}\left(  \mu\right)   &  \rightarrow\mathcal{O}_{\mu}\times
T^{\ast}A\label{E:defpsi}\\
\left(  g,\nu,a,\delta\right)   &  \longmapsto\left(  \nu,a,\delta\right)
\ \nonumber
\end{align}
is a surjective submersion that descends to diffeomorphism
\begin{align*}
\bar{\theta}:J_{G}^{-1}\left(  \mu\right)  /G_{\mu} &  \longrightarrow
\mathcal{O}_{\mu}\times T^{\ast}A\\
\left[  g,\nu,a,\delta\right]  _{G_{\mu}} &  \longmapsto\left(  \nu
,a,\delta\right) 
\end{align*}
that is symplectic with respect to the reduced symplectic form on the left and
\(\omega_{\mathcal{O}_{\mu}}^{-}+\omega_{T^{\ast}A}\) on the right.
The pushed-forward \(K\) action is \(k\cdot\left(  \nu,a,\delta\right)  =\left(
\Ad_{k^{-1}}^{\ast}\nu,k\cdot a,k\cdot\delta\right)  \).
It has momentum map
 \(J_{K}^{\prime}\left(  \nu
,a,\alpha\right)  =-\nu|_{\mathfrak{k}}+a\diamond_{\mathfrak{k}}\alpha.\)
\end{proposition}

Our results so far, combined with singular commution reduction (Theorem \ref{SCR}), 
give the following normal form for reduced spaces for cotangent-lifted actions. 

\begin{theorem}
\label{mainGK}Let \(G\) act properly on a manifold \(Q\) and by cotangent lifts on
\(T^{\ast}Q\) with momentum map \(J.\) Let \(z\in T_{q}^{\ast}Q\) and \(K=G_{q}.\) Let
\(\mu=J\left(  z\right)  \) and suppose that \(G_{\mu}\) is compact and
\(\mathcal{O}_{\mu}\) is locally closed. Let \(A\) be a \(K\)-invariant complement
to \(\mathfrak{g}\cdot q\) with respect to some \(K\)-invariant metric. Then there
is a local Poisson diffeomorphism between \(J^{-1}\left(  \mu\right)  /G_{\mu}\)
and the reduced space at \(0\) for the product action of \(K\) on the space
\(\mathcal{O}_{\mu}\times T^{\ast}A\) with symplectic form \(\omega
_{\mathcal{O}_{\mu}}^{-}+\omega_{T^{\ast}A},\) where \(K\) has the coadjoint
action on \(\mathcal{O}_{\mu}\) and the cotangent lifted action on \(T^{\ast}A.\)
\end{theorem}

Our main aim in this section is to characterise the symplectic normal space $N_s(z).$
To this end, we apply tangent-level commuting reduction (Theorem \ref{TLCR}) to the actions of $G$
and $K$ on $T^{\ast}\left(G\times A\right).$
Recall from Lemma \ref{xalpha} and Remark \ref{H} that $x=\left(e,\mu,0,\alpha\right)$ and
$H=G_z=G_{\left[x\right]_K}=G_{\mu}\cap K_{\alpha}.$ 
It is easy to check that \( \left(  G\times K\right)  _{x}=\{ (h,h) \, \big| \, h\in H\}, \)
and that $H=K_{\left[x\right]_{G_\mu}}$ as well.
A generalisation of this observation appears in Theorem \ref{TLCR}. As in that theorem, 
we will identify $\left(  G\times K\right)  _{x}$ with $H.$
The subgroup \(H\) acts on all three of the symplectic
normal spaces \(N_{s}\left(  \left[  x\right]  _{G_{\mu}}\right)  ,N_{s}\left(
x\right)  \) and \(N_{s}\left(  \left[  x\right]  _{K}\right)  \) in the usual
way. Theorem \ref{TLCR} (tangent-level commuting reduction) 
implies that the following maps are \(H\)-equivariant vector space
symplectomorphisms,
\begin{equation}
N_{s}\left(  \left[  x\right]  _{G_{\mu}}\right)  \overset{\overline{T_{x}\pi_{G_{\mu}}}}{\underset{\cong}{\longleftarrow}}N_{s}\left(  x\right)
\underset{\cong}{\overset{\overline{T_{x}\pi_{K}}}{\longrightarrow}}N_{s}\left(  \left[  x\right]  _{K}\right),  \label{E:TLCR0}\end{equation}
where the overbars denote the quotient maps, as in Remark \ref{fbar}.

Now recall from Equation \ref{E:phis} that
$T^{\ast}s^{-1}\circ\bar{\varphi}$ is a $G$-equivariant symplectomorphism from a neighbourhood
of $[x]_K$ to a neighourhood of $z.$
It follows that $\overline{T_{\left[x\right]_{K}}\left(T^{\ast}s^{-1}\circ\bar{\varphi}\right)}$
is an $H$-equivariant symplectomorphism from 
$N_{s}\left(\left[x\right]_{K}\right)$ to $N_{s}\left(z\right).$
Since $\varphi=\bar{\varphi}\circ\pi_{K},$ we can compose this
with $\overline{T_{x}\pi_{K}}$ from above to give 
$\overline{T_{x}\left(T^{\ast}s^{-1}\circ\varphi\right)}:N_s(x)\to N_s(z).$
Similarly, Proposition \ref{psi} implies that the map \(
\overline{T_{\left[x\right]_{G_{\mu}}}\bar{\theta}}:N_{s}\left(\left[x\right]_{G_{\mu}}\right)\rightarrow N_{s}\left(\mu,0,\alpha\right)\)
is an $H$-equivariant symplectomorphism; and we can compose this with $\overline{T_{x}\pi_{G_{\mu}}} $ to give 
 \(\overline{T_{x}\theta}:N_{s}\left(
x\right)  \rightarrow N_{s}\left(  z\right) . \)
Combining these results, we have the following,

\begin{theorem}
\label{mainVm}In the above context (with \(s, \varphi\) and \(\theta\) defined by Equations
\ref{E:defs}, \ref{E:cbphi2} and \ref{E:defpsi}, respectively), the composition \(\overline{T_{x}\theta}\circ\overline
{T_{x}\left(  T^{\ast}s^{-1}\circ\varphi\right)  }^{-1}:N_{s}\left(  z\right)
\rightarrow N_{s}\left(  \mu,0,\alpha\right)  \) is an \(H\)-invariant
symplectomorphism of symplectic normal spaces.
\end{theorem}

The space $N_{s}\left(\mu,0,\alpha\right)$ has
simple forms in the special cases $K\subset G_\mu$ and $\alpha=0.$
When $K\subset G_{\mu},$ the $K$ action on $\mathcal{O}_{\mu}$ is trivial, so
\begin{align}\label{E:KsubGmu}
N_{s}\left(\mu,0,\alpha\right)\cong T_\mu \mathcal{O}_\mu \oplus N_{s}\left(0,\alpha\right),
\end{align}
the second summand being 
the symplectic normal space at $\left(0, \alpha\right)$ for the cotangent-lifted action 
of $K$ on $T^{\ast}A.$ Recall that the momentum map 
for the latter action is
$J_A\left(a,\gamma\right)=a\diamond \gamma.$
It follows that $dJ_{A}\left(0,\alpha\right)\left(b,\beta\right)=0\diamond\beta+b\diamond\alpha
=b\diamond\alpha,$
so $\ker dJ_{A}\left(  0,\alpha\right)=\left(  \mathfrak{k}\cdot\alpha\right)  ^{\circ}\oplus A^*.$
Hence $N_{s}\left(0,\alpha\right)=
\ker dJ_A \left(  0,\alpha\right)  / \left(  \mathfrak{k}\cdot\left(
0,\alpha\right)  \right)  \cong\left(  \mathfrak{k}\cdot\alpha\right)  ^{\circ
}\oplus\left(  A^{\ast}/\left(  \mathfrak{k}\cdot\alpha\right)  \right).$
It is not hard to show that the dual \(\iota^{\ast
}\) of the inclusion \(\iota:\left(  \mathfrak{k}\cdot\alpha\right)  ^{\circ
}\hookrightarrow A\) descends to an isomorphism
$
\overline{\iota^{\ast}}:A^{\ast}/\left(  \mathfrak{k}\cdot\alpha\right)
\cong\left(  \left(  \mathfrak{k}\cdot\alpha\right)  ^{\circ}\right)  ^{\ast},
$
and the map
\begin{align}
 \label{E:linNs}
N_{s}\left(0,\alpha\right)\cong\left(  \mathfrak{k}\cdot\alpha\right)  ^{\circ
}\oplus\left(  A^{\ast}/\left(  \mathfrak{k}\cdot\alpha\right)  \right)
\overset{\left(  id,\overline{\iota^{\ast}}\right)  }{\longrightarrow}\left(
\mathfrak{k}\cdot\alpha\right)  ^{\circ}\oplus \left(  \left(  \mathfrak{k}\cdot\alpha\right)  ^{\circ}\right)  ^{\ast}\cong T^{\ast}\left(
\mathfrak{k}\cdot\alpha\right)  ^{\circ}
\end{align}
is a symplectomorphism and is $H$-equivariant with respect to the cotangent lift of the restriction of
the $K$ action on $A$ to a $H$ action on 
\(\left(  \mathfrak{k}\cdot\alpha\right)  ^{\circ}.\)
Thus we arrive at the following corollary to Theorem \ref{mainVm}:

\begin{corollary}
\label{NsKsubGmu}When \(K\subset G_{\mu},\) there is an \(H\)-equivariant
symplectomorphism
\[
N_{s}(z) \cong T_\mu\mathcal{O}_\mu \oplus T^{\ast}\left(
\mathfrak{k}\cdot\alpha\right)  ^{\circ}.
\]
%where \(T^{\ast}\left(  \mathfrak{k}\cdot\alpha\right)  ^{\circ}\) has the
%canonical cotangent bundle symplectic form and the \(H\) action is the cotangent
%lift of the \(K\) action on \(A\) restricted to an \(H\) action on \(\left(
%\mathfrak{k}\cdot\alpha\right)  ^{\circ}.\)
\end{corollary}

\begin{remark}
In light of Lemma \ref{cbiso} (iv), the above result applies whenever $K$ is normal in $G.$ 
\end{remark}

\begin{remark}
This corollary generalises a splitting established for free actions by
Montgomery et al. (see \cite{MMR84}).
\end{remark}%

We now consider the case $\alpha=0.$ Recall that $\alpha=\left. z\right|_A,$ where 
$A=\left(\mathfrak{g}\cdot q\right)^\perp; $ so, with respect to our choice of metric, this is the
case where the conjugate momentum $z$ is ``purely in the group direction.''
Since $J_{K}^{\prime}\left(\nu,a,\alpha\right)=-\nu|_{\mathfrak{k}}+a\diamond\alpha,$
it follows that $dJ_{K}^{\prime}\left(\mu,0,0\right)\left(\rho,b,\beta\right)=-\left.\rho\right|_{\mathfrak{k}}.$
Note that this equals $dJ_\mu\left(\mu\right)\left(\rho\right),$ where $J_\mu$ is the momentum
map for the coadjoint action of $K$ on $\left(\mathcal{O}_\mu,\omega_{\mathcal{O}_\mu}^-\right)$
namely $J_\mu(\nu)=-\left.\nu\right|_\mathfrak{k}.$
So $\ker dJ_{K}^{\prime}\left(\mu,0,0\right)=\ker dJ_{\mu}\left(\mu\right)\oplus T^*A.$
By Lemma \ref{cbiso} (iii), $J_{K}^{\prime}\left(\mu\right)=J_\mu(\mu) = 0,$
so $\mathfrak{k}_{J_{K}^{\prime}\left(\mu\right)} = \mathfrak{k}_{J_\mu(\mu)}
=\mathfrak{k}.$
It follows that 
\begin{align*}
N_s\left(\left(\mu,0,0\right)\right) 
&= \left(\ker dJ_{\mu}\left(\mu\right)\oplus T^*A\right) 
/ \left( \mathfrak{k}_{J_\mu(\mu)\cdot \mu} \oplus \{(0,0)\}\right)
= N_s(\mu) \oplus T^*A
\end{align*}
By the Reduction Lemma (or direct calculation), 
$N_s(\mu) = \left( \mathfrak{k}\cdot \mu \right)^{\omega^{-}} / \left( \mathfrak{k}\cdot \mu \right).$
Thus we have the following corollary to Theorem \ref{mainVm}:

\begin{corollary} \label{Nsalpha0}
When \(\alpha=0,\) the map in Theorem \ref{mainVm}
is an $H$-equivariant symplectomorphism 
\[
N_s(z)\cong N_s(\mu) \oplus T^{\ast}A
=\left( \mathfrak{k}\cdot \mu \right)^{\omega^{-}} / \left( \mathfrak{k}\cdot \mu \right)\oplus T^{\ast}A,
\]
where $N_s(\mu)$ is the symplectic normal space at $\mu$ for the 
coadjoint action of $K$ on $\left(\mathcal{O}_\mu,\omega_{\mathcal{O}_\mu}^-\right).$
\end{corollary}

\begin{remark} \label{releq}
The above corollary applies to all relative equilibria of simple mechanical systems.
Indeed, if $z\in T^*_qQ$ is such a relative equilbrium then $z=\mathbb{F}L\left( \xi \cdot q \right)$ 
for some $\xi\in\mathfrak{g}$
(see \cite{Mars92}). For any $v\in A$ we have 
$\left<z,v\right> = \left<\left< \xi\cdot q, v\right>\right>=0,$
since $A=\left(\mathfrak{g}\cdot q\right)^\perp.$
Hence $\alpha=\left. z \right|_A = 0.$
More generally, the corollary applies to any point $z$ such that the kernel of $z$ 
includes some complement to $\mathfrak{g}\cdot q,$ because we can
choose our metric on $Q$ such that this complement is $\left(\mathfrak{g}\cdot q\right)^\perp.$
\end{remark}

\begin{remark} \label{CBSLNFred} 
Theorem \ref{mainVm} and its corollaries, when combined with Theorem \ref{SLNFred},
give cotangent-bundle-specific local models of symplectic reduced spaces.
Indeed, if \(\mathcal{O}_{\mu}\) is
locally closed, there is a local Poisson diffeomorphism between the reduced
space \(J^{-1}\left(  \mu\right)  /G_{\mu}\) and the reduced space at \(0\) for
the \(H\) action on $N_s\left(\mu,0,\alpha\right).$
%If \(K\subset G_\mu\) 
%the latter space
%is \(T_\mu\mathcal{O}_\mu \oplus T^{\ast}\left(  \mathfrak{k}\cdot\alpha\right)  ^{\circ}\);
%if \(\alpha=0,\)
%it is $N_s(\mu) \oplus T^{\ast}A.$
Note the similarity to Theorem \ref{mainGK}, which shows that 
 \(J^{-1}\left(  \mu\right)  /G_{\mu}\) is isomorphic to the reduced space at \(0\) for
the \(K\) action on $\mathcal{O}_\mu\times T^*A.$
Thus, symplectic reduced spaces for cotangent bundles have two local models,
corresponding to the two isotropy subgroups $H$ and $K.$
The model involving $H$ and $N_s\left(\mu,0,\alpha\right)$ is more ``economical''
in that $H$ and $N_s\left(\mu,0,\alpha\right)$ may be smaller than 
$K$  and $\mathcal{O}_\mu\times T^*A,$ but on the other hand the latter space 
is ``simpler'' and might be easier to work with in some situations.
\end{remark}

We end this section with the observation that Corollaries \ref{NsKsubGmu}
and \ref{Nsalpha0} lead to refinements of the so-called 
 \emph{reconstruction equations} or  \emph{bundle equations} \cite{Ort98,OR02b,RWL02},
 which are a normal form for Hamilton's
equations in the coordinates given by the Hamiltonian
Slice Theorem (Theorem \ref{marle}).
Consider the local symplectomorphism 
$G\times_{H}\left(\mathfrak{m}^{\ast}\times N_{s}\right)\rightarrow P$
given by the Hamiltonian Slice Theorem,
for any proper globally Hamiltonian action of $G$ on $P,$
with $H=G_z$ as before;
recall that $\mathfrak{m}$ is an $H$-invariant complement to $\mathfrak{h}$
in $\mathfrak{g}_\mu,$ where $\mu=J(z).$
A Hamiltonian on $P$ pulls
back to a Hamiltonian $h$ on a neighbourhood of $\left[e,0,0\right]_{H}$
in $G\times_{H}\left(\mathfrak{m}^{\ast}\times N_{s}\right),$ with 
corresponding Hamiltonian vector field $X_{h}$. Using a local
bundle chart around $\left[e,0,0\right]_{H}$ for the principal bundle
$\pi_{H}:G\times\mathfrak{m}^{\ast}\times N_{s}\rightarrow G\times_{H}\left(\mathfrak{m}^{\ast}\times N_{s}\right),$
we can lift $X_{h}$ to a smooth vector field on a neighbourhood of
$\left(e,0,0\right)$ in $G\times\mathfrak{m}^{\ast}\times N_{s}.$
This lift is not unique;  however we can specify
a unique lift by 
choosing an $H$-invariant complement $\mathfrak{q}$ to $\mathfrak{g}_\mu,$
so that we now have $\mathfrak{g}=\mathfrak{q}+\mathfrak{m}+\mathfrak{h},$
and requiring that the component of the lifted vector field in the $\mathfrak{h}$ direction be zero.
The lifted vector field can now be written as
$X=\left(TL_{g}\left(X_{\mathfrak{m}}+X_{\mathrm{q}}\right),
X_{\mathfrak{m}^{\ast}},X_{N_{s}}\right).$

If it is possible to choose $\mathfrak{q}$ to be $G_\mu$-invariant then we say that $\mu$ is \emph{split}. 
We assume this now for simplicity; the general case is considered in \cite{RWL02}.
Assuming $\mu$ is split, it can be shown that the lifted Hamiltonian vector field $X(g,\rho,v)$
is given by:
\begin{align*}
X_{\mathfrak{q}} & =0\\
X_{\mathfrak{m}} & =D_{\mathfrak{m}^{\ast}}\left(h\circ\pi\right)(\rho,v)\\
\mathrm{i}_{X_{N_{s}}}\omega_{N_{s}}&=D_{N_{s}}\left(h\circ\pi\right)(\rho,v)\\
X_{\mathfrak{m}^{\ast}} & =\mathbb{P}_{\mathfrak{m}^{\ast}}\left(\ad _{D_{\mathfrak{m}^{\ast}}\left(h\circ\pi\right)}^{\ast}\rho\right)+\ad _{D_{\mathfrak{m}^{\ast}}\left(h\circ\pi\right)}^{\ast}J_{N_{s}}\left(v\right)\end{align*}

Now suppose $P=T^*Q$ and $z\in T_q^*Q.$ 
If $G_q\subset G_{\mu},$ we know from Corollary \ref{NsKsubGmu} that
$N_{s}$ is linearly symplectomorphic to 
$T_\mu\mathcal{O}_\mu\oplus T^{\ast}\left(\mathfrak{k}\cdot\alpha\right)^{\circ}.$
Let $B=\left(\mathfrak{k}\cdot\alpha\right)^{\circ},$
so $N_s\cong T_\mu\mathcal{O}_\mu\oplus B\oplus B^*.$
The vector field $X_{N_{s}}$ separates into three components $X_\mu,X_{B}$
and $X_{B^{\ast}}$ and the equation 
$\mathrm{i}_{X_{N_{s}}}\omega_{N_{s}}=D_{N_{s}}\left(h\circ\pi\right)(\rho,v)$
splits into the three equations, 
\begin{align} \label{E:reconGmuG}
\mathrm{i}_{X_{\mu}}\omega_{\mathcal{O}_\mu^-}
&=D_{T\mu\mathcal{O}_\mu}\left(h\circ\pi\right)(\rho,v)\\
X_{B} & =D_{B^{\ast}}\left(h\circ\pi\right)(\rho,v) \notag \\
X_{B^{\ast}} & =-D_{B}\left(h\circ\pi\right)(\rho,v). \notag
\end{align}
(the last two equations being the canonical form for Hamilton's equations).
The case $\alpha=0$ in Corollary \ref{Nsalpha0} is similar: the $X_{N_s}$ equation
is replaced by three equations,
\begin{align} \label{E:reconalpha0}
\mathrm{i}_{X_{\mu}}\omega_{\mathcal{O}_\mu^-,red}&
=D_{N_{s}(\mu)}\left(h\circ\pi\right)(\rho,v)\\
X_{A} & =D_{A^{\ast}}\left(h\circ\pi\right)(\rho,v) \notag \\
X_{A^{\ast}} & =-D_{A}\left(h\circ\pi\right)(\rho,v). \notag
\end{align}

\section{A cotangent bundle slice theorem\label{sectcbs}}

In this section, we extend the Hamiltonian slice theorem (Theorem \ref{marle})
in the case of a lifted action on a cotangent bundle.
The main result is Theorem \ref{mainGmuG}.
We will consider only the case of fully isotropic momenta, $G_{\mu}=G,$ for reasons that
will be summarised in Remark \ref{whyGmuG}.
Our model for $T^*Q$ will be $G\times_H\left(\mathfrak{m}^* \oplus N_s\right),$
as in the general Hamiltonian theorem,
with the same symplectic form as in that theorem
(definitions will be reiterated below).
However, in contrast to the general Hamiltonian slice theorem,
our isomorphism from the model space to $T^*Q$
will be constructed explicitly, apart from the use of a Riemannian exponential in the base space $Q$.
The construction will use the decomposition of $N_s$ in Corollary \ref{NsKsubGmu}.

As before, let $G$ be a Lie group acting smoothly and properly by cotangent lifts 
on $T^*Q,$ with momentum map $J.$ Let $z\in T^*_qQ$ and $\mu=J(z),$ and let $K=G_{q}$
and $H=G_z.$ We assume $G_{\mu}=G.$
Let $\mathfrak{g},\mathfrak{g}_\mu,\mathfrak{h}$ and $\mathfrak{k}$ be the Lie algebras
of $G,G_{\mu},H$ and $K$ respectively. 
Fix a $K$-invariant inner product on $\mathfrak{g}$ and let 
$\mathfrak{m}=\mathfrak{k}^\perp.$ Let $N_{s}$ be the
symplectic normal space at $z.$ Our goal is to find a symplectic tube from 
$G\times_H\left(\mathfrak{m}^* \oplus N_s\right)$ to $T^*Q$ that takes $[e,0,0]$ to $z.$
%where the domain has the same symplectic form as in the Hamiltonian slice theorem.

We first apply Palais' slice theorem in the configuration space $Q.$
Fix a $K$-invariant Riemannian metric on $Q$ and let $A=\left(\mathfrak{g}\cdot q\right)^\perp.$
By the slice theorem, there exists a
$G$-equivariant diffeomorphism $s:G\times_{K}V\to Q$ taking $[e,0]_K$ to $q,$
for some neighbourhood $V$ of $0$ in $A.$
The cotangent lift of $s$ is a
$G$-equivariant symplectomorphism
$T^*s^{-1}:T^{\ast}\left(G\times_{K}V\right)\to T^*Q.$
Let $\varphi:J_{K}^{-1}\left(0\right)\to T^{\ast}\left(G\times_{K}A\right)$ be the cotangent bundle reduction map defined in Equation \ref{E:cbphi2},
 and let $\alpha=\left. z\right|_A.$
Recall from Lemma \ref{xalpha} that $T^*s^{-1}\circ \varphi\left(e,\mu,0,\alpha\right)=z.$
Hence it will suffice to find a symplectic tube 
%(i.e., a $G$-equivariant local symplectomorphism)
\begin{align}\label{E:deftau}
\tau : G\times_{H}\left(\mathfrak{m}^*\oplus N_s\right) & \longrightarrow T^{\ast}\left(G\times_{K}A\right)\\
\left[e,0,0\right]_H & \longmapsto\varphi\left(e,\mu,0,\alpha\right).\nonumber \end{align}

Since $G_{\mu}=G,$ Corollary \ref{NsKsubGmu} says that the symplectic
normal space $N_{s}$ is $H$-equivariantly symplectomorphic
to $T^{\ast}B\cong B\oplus B^{\ast},$ where $B=\left(\mathfrak{k}\cdot\alpha\right)^{\circ}\subset A;$
the symplectic form on $T^{\ast}B$ is the canonical one, and the $H$ action on $T^{\ast}B$ is the
cotangent-lift of the restriction to $H$ and $B$ of the $K$ action on $A$.
We will identify $N_s$ with $B\oplus B^*.$
Recall that the presymplectic form on $G\times_{H}\left(\mathfrak{m}^*\oplus N_s\right)$ 
in the Hamiltonian slice theorem is defined using a symplectic form 
$\Omega_{Z}$ $=\Omega_{c}+\Omega_{\mu}+\Omega_{N_{s}}$ on $Z=G\times\mathfrak{g}_{\mu}^{\ast}\times N_{s}$ (see Equation \ref{E:L}).
Since we are assuming $G_\mu=G,$ we have
$Z=G\times\mathfrak{g}^{\ast}\times B\times B^{\ast},$
which we identify with $T^{\ast}\left(G\times B\right)$ by left trivialisation
of $T^{\ast}G.$ 
The twist action of $H$ on $Z$ becomes the cotangent lift of the twist action of $H$ on
$G\times B.$
The form $\Omega_{\mu}$ is a pull-back of a symplectic form on $\mathcal{O}_{\mu},$ which
is trivial in this case, so $\Omega_{Z}=\Omega_{c}+\Omega_{N_{s}}.$ 
Since $\Omega_{c}$ is the pull-back by
left-trivialisation of the canonical symplectic form on $T^{\ast}G,$ and $\Omega_{N_{s}}$
is the canonical symplectic form on $T^{\ast}B,$ the identification
of $Z$ with $T^{\ast}\left(G\times B\right)$ makes $\Omega_{Z}$
the canonical symplectic form on $T^{\ast}\left(G\times B\right).$
Note that, unlike in the general case, this $\Omega_Z$ is nondegenerate everywhere.

The symplectic form on $G\times_{H}\left(\mathfrak{m}^*\oplus N_s\right)$ is defined via 
an isomorphism with $J_H^{-1}(0)/H,$ where $J_H$ is the momentum map of the $H$ action
on $Z=T^{\ast}\left(G\times B\right).$ The isomorphism, defined earlier in Equations \ref{E:l} and \ref{E:L}, is
\begin{align*}
L:G\times_H\left(\mathfrak{m}^{\ast}\oplus B \oplus B^*\right)
& \longrightarrow J_{H}^{-1}\left(0\right)/H\\
\left[g,\nu,a,\delta\right]_H & \longmapsto
\left[g,\nu + a \diamond_\mathfrak{h} \delta,a,\delta\right]_H\nonumber \end{align*}
The symplectic form on $G\times_H\left(\mathfrak{m}^{\ast}\oplus B \oplus B^*\right)$
is defined as the pull-back by $L$ of the reduced symplectic form on $J_{H}^{-1}\left(0\right)/H.$
Since $L$ is clearly $G$-equivariant, it is a symplectic tube.

In the present case, cotangent bundle reduction (Theorem \ref{RPCBRZ}) 
shows that $J_H^{-1}(0)/H$ is isomorphic to $T^{\ast}\left(G\times_H B\right).$
Let  $\psi$ and $\bar{\psi}$ be the maps in the cotangent bundle reduction theorem,
\begin{equation*}
\begin{array}{ccccc}
%G\times\left(\mathfrak{m}^{\ast}\oplus B\oplus B^{\ast}\right) & \overset{l}{\longrightarrow} & 
J_{H}^{-1}\left(0\right) & \hookrightarrow & T^{\ast}\left(G\times B\right)\\
%\downarrow\pi_{H} &  &
 \pi_{Z,H}\downarrow & \overset{\psi}{\searrow}\\
%G\times_{H}\left(\mathfrak{m}^{\ast}\oplus B\oplus B^{\ast}\right) & \overset{L}{\longrightarrow} & 
J_{H}^{-1}\left(0\right)/H & \overset{\bar{\psi}}{\longrightarrow} & T^{\ast}\left(G\times_{H}B\right)\end{array}\end{equation*}
It is easily checked that $\bar\psi$ is $G$-equivariant, by the same reasoning as used in Proposition \ref{unrollK}.
%In summary, if we define the map $l$ as in Equation \ref{E:l}, namely
Note that
\begin{equation}
\left(\bar{\psi}\circ
L\right)\left(\left[g,\nu,a,\delta\right]_{H}\right)=\psi\left(g,\nu+a\diamond_{\mathfrak{h}}\delta,a,\delta
\right).\label{E:Lpsi}\end{equation}
In particular, $\left(\bar{\psi}\circ L\right)\left(\left[e,0,0,0\right]_{H}\right)=\psi\left(e,0,0,0\right).$
Thus, to find a tube $\tau$ as in Equation \ref{E:deftau}, it suffices to find
a $G$-equivariant symplectomorphism
\begin{align}
\label{E:defsigmabar}
\bar{\sigma}:T^{\ast}\left(G\times_{H}B\right) & \longrightarrow T^{\ast}\left(G\times_{K}A\right)\\
\psi\left(e,0,0,0\right) & \longmapsto\varphi\left(e,\mu,0,\alpha\right)\nonumber \end{align}
(we will have to restrict the domain of $\bar\sigma$ in the general case).

\noindent \textbf{The $z=0$ case.} 
In the simplest case, $z=0\in T^*_qQ,$
we have $\mu=\alpha = 0,\ H=K,\, B=A$ and $\psi=\varphi$, so $\bar{\sigma}$ may be
chosen to be the identity map on $T^*(G\times_KA).$
Composing this with the maps $\bar{\psi}\circ L$ and $T^*s^{-1}$ gives 
the symplectic tube 
\begin{align*}
G\times_{K}\left(\mathfrak{m}^{\ast}\oplus V\oplus A^{\ast}\right)
& \longrightarrow T^*Q \\
\left[g,\nu,a,\delta\right]_K &\longmapsto T^*s^{-1} \circ \varphi 
\left(g,\nu+a\diamond_{\mathfrak{h}}\delta,a,\delta\right), \notag
\end{align*}
where  $V$ is the neighbourhood of $0$ in $A$ given by Palais' slice theorem applied at $q\in Q.$
%This tube has been known for some time in the context of
%molecular chemistry (see \cite{Wat70,KRT00}).

\medskip
In general, $\mu$ and $\alpha$ may be nonzero. However, if $H=K,$ we may take
$\bar\sigma$ to be a simple shift map, as in the following lemma. Note that,
since $H=K_\alpha$ and $B=\left(\mathfrak{k}\cdot \alpha\right)^\circ,$
the condition $H=K$ is equivalent to $B=A.$

\begin{lemma}
\label{shift} If $G_{\mu} = G$ then the shift map $\Sigma_{(\mu,\alpha)}:
\left(  g,\nu,a,\delta\right)   \mapsto\left(  g,\mu+\nu,a,\alpha+\delta
\right),$
 from $T^*(G\times A)$ to itself, is symplectic
 and $G$-equivariant.
If $H=K$ then $B=A$ and $\Sigma_{(\mu,\alpha)}$
leaves \(J_{H}\) invariant and is $H$-equivariant.
The ``quotient'' of $\Sigma_{(\mu,\alpha)}$ by $\psi,$ the map
\begin{align*}
\overline{\Sigma}_{(\mu,\alpha)}:
T^*\left(G\times_H B\right) &\longrightarrow T^*\left(G\times_H B\right) \\
\psi \left(g,\nu,a,\delta\right)  &\longmapsto \psi \left(g,\mu+\nu,a, \alpha + \delta\right)
\end{align*}
is a $G$-equivariant symplectomorphism.
\end{lemma}

\begin{proof}
It is clear from the local coordinate formula $dq^i\wedge dp_i$ that a canonical cotangent bundle
symplectic form is invariant under shifts in the $p$ variable.
The $G$-equivariance is also clear.  
Now suppose $H=K,$ which implies $B=A,$ as explained above.
Recall that $J_H(g,\nu,a,\delta)=-\left.\nu\right|_\mathfrak{h} + a\diamond_\mathfrak{h} \delta.$
Since \(\mathfrak{k}\subset\ker\mu\) (see Lemma \ref{cbiso} (iii)),
it follows that  $-\left.\mu \right|_\mathfrak{h}= 0$.
From $H= K_\alpha$  (see Remark \ref{H}) it follows that 
$a\diamond_\mathfrak{h} \alpha=0.$ Hence $J_H$ is invariant under $\Sigma_{(\mu,\alpha)}.$
The $H$-equivariance follows from the linearity of the $H$
actions on $\mathfrak{g}^*$ and $B^*$ together with $H\subset G_\mu\cap K_\alpha.$

Since $\Sigma_{(\mu,\alpha)}$ is a 
$G$- and  $H$- equivariant symplectomorphism leaving $J_H^{-1}(0)$ invariant, 
it descends to a $G$-equivariant symplectomorphism from $J_H^{-1}(0)/H$ to itself.
This map induces $\overline{\Sigma}_{(\mu,\alpha)}$ via the 
$G$-equivariant symplectomorphism $\overline{\varphi}:J_H^{-1}(0)/H\to T^*(G\times_H B).$
\end{proof}

%\medskip
%This shift map is sufficient to 
%solve our problem whenever $B=A.$ 
%(Recall that $H=G_z,$  the phase space isotropy group,
%while $K=G_q,$ the configuration space isotropy  group.)

\medskip

\noindent \textbf{The case }$\mathbf{H=K}.$
Subcases include: $z=0;$ $\alpha=0;$ and all relative equilibria of simple mechanical systems (see Remark \ref{releq})
(recall that we are assuming $G_\mu = G$ throughout this section).
When $H=K$ we have $B= A,$ so 
$T^*(G\times_HB)=T^*(G\times_KA)$ and $\psi = \varphi$. Thus we may take
$\bar{\sigma}=\overline{\Sigma}_{(\mu,\alpha)}.$
Composing this with the maps $\bar{\psi}\circ L$ and $T^*s^{-1}$ gives
\begin{align*}
G\times_{K}\left(\mathfrak{m}^{\ast}\oplus V\oplus A^{\ast}\right)
& \longrightarrow T^*Q \\
\left[g,\nu,a,\delta\right]_K &\longmapsto T^*s^{-1} \circ \varphi 
\left(g,\mu+\nu+b\diamond_{\mathfrak{h}}\delta,b,\alpha +\delta\right),
\end{align*}
where  $V$ is the neighbourhood of $0$ in $A$ given by Palais' slice theorem applied at $q\in Q.$

\medskip

The general case is more difficult.
We identify $B^*$ with $\left(\mathfrak{k}\cdot \alpha\right)^\perp\subset A^*.$
It is easily checked that the shift formula 
$\left(  g,\nu,a,\delta\right)   \mapsto\left(  g,\mu+\nu,a,\alpha+\delta
\right),$
 as a map from $T^*(G\times B)$ to $T^*(G\times A),$
 need \emph{not} map $J_H^{-1}(0)$ into $J_K^{-1}(0),$
 so cannot be used directly to define a map $\bar{\sigma}$ as in Equation \ref{E:defsigmabar}.
%It is tempting to try the cotangent lift of the map $\left[g,b\right]_{H}\rightarrow\left[g,b\right]_{K}.$
%However, the cotangent lift of this or \emph{any} map taking $\left[e,0\right]_{H}$
%$\,\,$to $\left[e,0\right]_{K}$ would map $\psi\left(e,0,0,0\right)$
%to $\varphi\left(e,0,0,0\right),$ not to $\varphi\left(e,\mu,0,\alpha\right)$
%as desired. We will see later in this section that the idea of
%using cotangent lifts can be made to work (see Theorem \ref{alttube}), but for now we follow a
%different approach: 
We will look for a map as close a possible to this shift map but
with image contained in $J_K^{-1}(0).$
We will conclude in Lemma \ref{sigma} that there is a unique map of the form
$\left(g,\nu,a,\delta\right)\mapsto\left(g,\mu+\nu,a+c,\alpha+\delta\right),$
for $c\in B^{\perp},$ that accomplishes this.

%as in the following proposition. 
%\begin{proposition}
%\label{dropsigma}Let \(U\) be an open subset of \(J_{H}^{-1}\left(  0\right)
%\subset Z=T^{\ast}\left(  G\times B\right)  \) and let \(\sigma:U\rightarrow
%\left(  J_{K}^{-1}\left(  0\right)  \subset T^{\ast}\left(  G\times A\right)
%\right)  \) be a \(G\)-equivariant presymplectic diffeomorphism with respect to
%the canonical symplectic forms on \(T^{\ast}\left(  G\times B\right)  \) and
%\(T^{\ast}\left(  G\times A\right)  .\) Then the following commuting diagram
%defines a unique \(G\)-equivariant symplectomorphism \(\tau_{1},\)
%\[\begin{array}
%[c]{ccc}\left(  U\subset J_{H}^{-1}\left(  0\right)  \right)  & \overset{\sigma
%}{\hookrightarrow} & J_{K}^{-1}\left(  0\right) \\
%\downarrow\psi &  & \downarrow\varphi\\
%\left(  \psi\left(  U\right)  \subset T^{\ast}\left(  G\times_{H}B\right)
%\right)  & \overset{\tau_{1}}{\hookrightarrow} & T^{\ast}\left(  G\times
%_{K}A\right)  .
%\end{array}
%\]
%\end{proposition}

We proceed by characterising the space $\left(G\times\mathfrak{g}^{\ast}\times A\times\left(\alpha+B^{\ast}\right)\right)\cap J_{K}^{-1}\left(0\right).$
We decompose $A$ as $\left(\left(\mathfrak{k}\cdot\alpha\right)^{\bot}\right)^{\circ}\oplus\left(\mathfrak{k}\cdot\alpha\right)^{\circ}=B^{\bot}\oplus B;$ this splitting is $H$-invariant, since \(H\) fixes \(\alpha\) (see Remark \ref{H}).
By definition of $\mathfrak{m},$ we have an $H$-equivariant splitting
$\mathfrak{g=m}\oplus\mathfrak{h}$ (recall that we are assuming $G_{\mu}=G$). 
It is easily checked that 
\(\mathfrak{k}\) splits \(H\)-equivariantly as \(\mathfrak{k}=\left(
\mathfrak{m}\cap\mathfrak{k}\right)  \oplus\mathfrak{h}\).

%\begin{proof}
%It is clear that \(\left(  \mathfrak{m}\cap\mathfrak{k}\right)  \) and \(\mathfrak{h}\)
%have trivial intersection. Now let \(\xi\in\mathfrak{k.}\) Since \(\mathfrak{g=m}\oplus\mathfrak{h,}\) we %can write \(\xi=\xi_{\mathfrak{m}}+\xi_{\mathfrak{h}},\) for
%some \(\xi_{\mathfrak{m}}\in\mathfrak{m}\) and \(\xi_{\mathfrak{h}}\in\mathfrak{h.}\)
%Since \(\mathfrak{h}\subset\mathfrak{k,}\) we have \(\xi_{\mathfrak{m}}=\xi
%-\xi_{\mathfrak{h}}\in\mathfrak{k,}\) so \(\xi_{\mathfrak{m}}\in\mathfrak{m}\cap\mathfrak{k.}\) The \(H\)-%equivariance follows from the \(H\)-equivariance of the
%original splitting \(\mathfrak{g=m}\oplus\mathfrak{h.}\)
%\end{proof}

\begin{lemma}
\label{imsigma}
\begin{align*}
&  \left(  G\times\mathfrak{g}^{\ast}\times A\times\left(  \alpha+B^{\ast
}\right)  \right)  \cap J_{K}^{-1}\left(  0\right) \\
&  =\left\{
\left(  g,\nu,b+c,\alpha+\delta\right)  \mid b\in B, c\in B^{\perp}, J_{H}\left(
g,\nu,b,\delta\right)  =0
\text{ and }
\left. \left( -\nu+b\diamond \delta+c\diamond\left(
\alpha+\delta\right) \right) \right|  _{\mathfrak{m}\cap\mathfrak{k}} =0
\right\}
\end{align*}
and this is a submanifold of $T^*(G\times A).$
\end{lemma}

\begin{proof}
Recall that \(J_{K}\left(  g,\nu,a,\beta\right)
=-\left.  \nu\right|  _{\mathfrak{k}}+a\diamond_{\mathfrak{k}}\beta\).
The restriction of $J_K$ to
$G\times\mathfrak{g}^{\ast}\times A\times\left(  \alpha+B^{\ast}\right)$ is a submersion,
since $\nu\mapsto \left.\nu\right|_\mathfrak{k}$ is one.
It follows that 
$ \left(  G\times\mathfrak{g}^{\ast}\times A\times\left(  \alpha+B^{\ast
}\right)  \right)  \cap J_{K}^{-1}\left(  0\right)$
is a submanifold of  $\left(  G\times\mathfrak{g}^{\ast}\times A\times\left(  \alpha+B^{\ast
}\right)  \right),$ and hence of $G\times\mathfrak{g}^{\ast}\times A\times A^*.$

Now let \(\left(  g,\nu,b+c,\alpha+\delta\right)  \in\left(
G\times\mathfrak{g}^{\ast}\oplus A\oplus\left(  \alpha+B^{\ast}\right)  \right)
,\) with \(b\in B\) and \(c\in B^{\perp}.\) 
Then \(J_{K}\left(  g,\nu,b+c,\alpha+\delta\right)
=-\left.  \nu\right|  _{\mathfrak{k}}+(b+c)\diamond_{\mathfrak{k}}(\alpha+\delta)\).
Since \(b\in B=\left(  \mathfrak{k}\cdot\alpha\right)  ^{\circ},\) it follows
that \(b\diamond_{\mathfrak{k}}\alpha=0.\)
For any \(\xi\in\mathfrak{h}\) 
we have \(\xi
\cdot\alpha=0\) and \(\xi\cdot\delta\in\left(  \mathfrak{k}\cdot\alpha\right)
^{\perp},\) and so \(\left\langle c,\xi\cdot\left(  \alpha+\delta\right)
\right\rangle =0;\) it follows that \(c\diamond_{\mathfrak{h}}\left(
\alpha+\delta\right)  =0.\) Thus
\begin{align*}
\left(  g,\nu,b+c,\alpha+\delta\right)   \in J_{K}^{-1}\left(  0\right) 
&  \Longleftrightarrow\left.  \nu\right|  _{\mathfrak{k}}=b\diamond_{\mathfrak{k}}\delta+c\diamond_{\mathfrak{m}\cap\mathfrak{k}}\left(  \alpha+\delta\right) \\
&  \Longleftrightarrow\left.  \nu\right|  _{\mathfrak{h}}-b\diamond_{\mathfrak{h}}\delta=-\left.  \nu\right|  _{\mathfrak{m}\cap\mathfrak{k}}+b\diamond
_{\mathfrak{m}\cap\mathfrak{k}}\delta+c\diamond_{\mathfrak{m}\cap\mathfrak{k}}\left(
\alpha+\delta\right) \\
&  \Longleftrightarrow\left.  \nu\right|  _{\mathfrak{h}}-b\diamond_{\mathfrak{h}}\delta=0\text{ and }-\left.  \nu\right|  _{\mathfrak{m}\cap\mathfrak{k}}+b\diamond_{\mathfrak{m}\cap\mathfrak{k}}\delta+c\diamond_{\mathfrak{m}\cap\mathfrak{k}}\left(  \alpha+\delta\right)  =0\\
&  \Longleftrightarrow J_{H}\left(  g,\nu,b,\delta\right)  =0\text{ and}
\left. \left( -\nu+b\diamond \delta+c\diamond\left(
\alpha+\delta\right) \right) \right|  _{\mathfrak{m}\cap\mathfrak{k}} =0.
\end{align*}
\end{proof}

%\marginpar{$\Gamma$ is like the mechanical connection the action of $K$ on
%$A^{\ast}$%
%}

\begin{lemma}
\label{sigma}Let \(U\) be an \(H\)-invariant neighbourhood of \(0\) in \(B^{\ast}\)
such that the map
\begin{align*}
t:K\times_{H}U  &  \longrightarrow A^{\ast}\\
\left[  k,\delta\right]  _{H}  &  \longmapsto k\cdot\left(\alpha + \delta\right)
\end{align*}
is injective; such a \(U\) always exists. Then

\begin{enumerate}
\item  For every \(\delta\in U,\) the map \(\Gamma_{\delta}^{\ast}:\left(
\mathfrak{m}\cap\mathfrak{k}\right)  ^{\ast}\longmapsto B^{\perp}\) defined by
\begin{equation}
\left\langle \Gamma_{\delta}^{\ast}\left(  \nu\right)  ,\xi\cdot\left(
\alpha+\delta\right)  +\varepsilon\right\rangle =\left\langle \nu
,\xi\right\rangle , \label{E:Gamma}\end{equation}
for every \(\xi\in\mathfrak{m}\cap\mathfrak{k}\) and \(\varepsilon\in B^*,\) 
is $H$-equivariant and has an inverse given by
\[
c\longmapsto-c\diamond_{\mathfrak{m}\cap\mathfrak{k}}\left(  \alpha+\delta\right)
.
\]

\item The map $\sigma$ defined by 
\begin{align*}
\label{E:sigma}
\sigma: & \left(G\times\mathfrak{g}^{\ast}\times B\times U\right)\cap J_{H}^{-1}\left(0\right)\longrightarrow\left(G\times\mathfrak{g}^{\ast}\times A\times\left(\alpha+U\right)\right)\cap J_{K}^{-1}\left(0\right)\\
\left(g,\nu,a,\delta\right) & \longmapsto\left(g,\mu+\nu,a+\Gamma_{\delta}^{\ast}\left(-\left.\nu\right|_{\mathfrak{m}\cap\mathfrak{k}}+a\diamond_{\mathfrak{m}\cap\mathfrak{k}}\delta\right),\alpha+\delta\right)\end{align*}
 is the unique function, with this domain and range, of the form $\left(g,\nu,a,\delta\right)\longmapsto\left(g,\mu+\nu,a+c,\alpha+\delta\right)$
for $c\in B^{\perp}.$ It is a presymplectic diffeomorphism,
with respect to the canonical symplectic forms  on
$T^{\ast}\left(G\times B\right)$ and $T^{\ast}\left(G\times A\right).$
It is equivariant with respect to both the left multiplication action by $G$ and
the twist action by $H.$
It descends to a \(G\)-equivariant symplectic embedding
\(\bar{\sigma}\) defined by the following commutative diagram,
where $W=\left(G\times\mathfrak{g}^{\ast}\times B\times U\right)\cap J_{H}^{-1}\left(0\right),$
\[\begin{array}
[c]{ccc}\left(  W\subset J_{H}^{-1}\left(  0\right)  \right)  & \overset{\sigma
}{\hookrightarrow} & J_{K}^{-1}\left(  0\right) \\
\downarrow\psi &  & \downarrow\varphi\\
\left(  \psi\left(  W\right)  \subset T^{\ast}\left(  G\times_{H}B\right)
\right)  & \overset{\bar{\sigma}}{\hookrightarrow} & T^{\ast}\left(  G\times
_{K}A\right) \ .
\end{array}
\]
The image of $\bar\sigma$ is an open subset of $T^{\ast}\left(  G\times_{K}A\right).$

\end{enumerate}
\end{lemma}

\begin{proof}
(1) Note that  \(H=K_{\alpha}\) and 
we are identifying \(B^{\ast}$ with $\left(  \mathfrak{k}\cdot\alpha\right)  ^{\bot},\)
which is a linear slice for the the \(K\) action on \(A^{\ast}.\) 
The slice theorem for linear actions (Theorem \ref{linslice}) shows the existence of an \(H\)-invariant neighbourhood \(U\) of \ \(0\) in
\(B^{\ast}\) such that the map $t$ above is injective and that, given any such $U$
the map \(t\) is a \(K\)-equivariant
diffeomorphism onto a \(K\)-invariant neighbourhood of \(\alpha.\) 
Let \(\pi
_{H}:K\times U\rightarrow K\times_{H}U\) be projection. The composition
\(t\circ\pi_{H}\) is a submersion. For any \(k\in K\) and \(\delta\in U,\) the
kernel of \(T_{\left(  k,\delta\right)  }\pi_{H}\) is \(\left\{  \left(
-\zeta,\zeta\cdot\delta\right)  \in\mathfrak{k}\oplus B^*
\mid\zeta\in\mathfrak{h}\right\}  ,\) which is a
complement to the space \(\left(  \mathfrak{m}\cap\mathfrak{k}\right)
\oplus B^*\) in \(\mathfrak{k}\oplus B^*,\) so
\begin{align*}
\left.  T_{\left(  e,\delta\right)  }\left(  t\circ\pi_{H}\right)  \right|
_{\left(  \mathfrak{m}\cap\mathfrak{k}\right)  \oplus B^*}
:\left(  \mathfrak{m}\cap\mathfrak{k}\right)  \oplus B^* 
&  \longrightarrow T_{\left(
\alpha+\delta\right)  }A^{\ast} \cong A^*\\
\left(  \xi,\varepsilon\right)   
&  \longmapsto  \xi\cdot\left( \alpha+\delta\right)  +\varepsilon
\end{align*}
is an isomorphism. 
It follows that Equation \ref{E:Gamma} defines a map \(\Gamma_{\delta}^{\ast}\)
from \(\left(  \mathfrak{m}\cap\mathfrak{k}\right)  ^{\ast}\) to \(A;\) its image
is clearly contained in \(\left(  B^*\right)  ^{\circ}\cong B^{\perp}.\)

It is easily checked that  \(\Gamma_{\delta}^{\ast}\) is
$H$-equivariant and has the stated inverse.
%f:B^{\perp}  &  \rightarrow\left(  \mathfrak{m}\cap\mathfrak{k}\right)  ^{\ast}\\
%$c   \longmapsto-c\diamond_{\mathfrak{m}\cap\mathfrak{k}}\left(  \alpha
%+\delta\right) $

(2) 
We first check that \(\sigma\) is well-defined; the only part that needs
checking is that its image is contained in the target space. 
It follows from Claim 1 that the condition  
$\left. \left( -\nu+b\diamond \delta+c\diamond\left(
\alpha+\delta\right) \right) \right|  _{\mathfrak{m}\cap\mathfrak{k}} =0$
in  Lemma \ref{imsigma}
is equivalent to 
$c=\Gamma_\delta^* \left(\left. \left( -\nu+b\diamond \delta\right)\right|_{\mathfrak{m}\cap\mathfrak{k}}
\right).$
The other condition in Lemma \ref{imsigma} that needs checking is
\(J_{H}\left(  g,\mu+\nu,a,\delta\right)  =0,\) for
every \(\left(  g,\nu,a,\delta\right)  \in\left(  G\times\mathfrak{g}^{\ast
}\times B\times U\right)  \cap J_{H}^{-1}\left(  0\right); \) 
this follows easily from the the fact that $\left.\mu\right|_\mathfrak{h} = 0.$
It is easily checked that $\sigma$ has an inverse given by 
$\left(  g,\mu+\nu,a+c,\alpha+\delta\right)  \mapsto\left(  g,\nu
,a,\delta\right),$
where \(a\in B\) and \(c\in B^{\perp}.\)

Part of that same argument, namely the fact that \(\Gamma_{\delta}^{\ast
}\left(  \left.  \nu\right|  _{\mathfrak{m}\cap\mathfrak{k}}-a\diamond
_{\mathfrak{m}\cap\mathfrak{k}}\delta\right)  =-c\) for any \(\left(  g,\mu
+\nu,a+c,\alpha+\delta\right)  \) in the range of \(\sigma\), also proves that
\(\sigma\) is the unique function, with the given domain and range, of the form
\(\left(  g,\nu,a,\delta\right)  \longmapsto\left(  g,\mu+\nu,a+c,\alpha
+\delta\right)  \) for \(c\in B^{\perp}.\)

We now show that \(\sigma\) is a diffeomorphism. Note  that
its domain is a submanifold of \(T^{\ast}\left(  G\times B\right)  ,\) being an
open subset of a level set of the momentum map of a free action;
similarly the range of \(\sigma\) is a submanifold of \(T^{\ast}\left(  G\times
A\right)  .\) 
Since the image of \(\Gamma_{\delta}^{\ast}\) is \(B^{\perp},\) its
derivative is always in \(B^{\perp},\) so for any \(\left(  g,\nu,a,\delta\right)  \) in the
domain of \(\sigma,\) and any tangent vector \(\left(  \dot{g},\dot{\nu},\dot
{a},\dot{\delta}\right)  \in\mathfrak{g}\times\mathfrak{g}^{\ast}\times B\times
B^{\ast},\) we have
$T_{\left(  g,\nu,a,\delta\right)  }\sigma\left(  \dot{g},\dot{\nu},\dot
{a},\dot{\delta}\right)  =\left(  \dot{g},\dot{\nu},\dot{a}+\dot{c},\dot{\delta}\right)$
%\quad \textrm{for some} \dot{c}\in B^{\perp}. 
for some $\dot{c}\in B^{\perp}.$ 
It is clear from this formula that  \(\sigma\) is an immersion. 
But any bijective immersion is a diffeomorphism (see Lemma \ref{bijimm} below).

We next show that \(\sigma\)
is presymplectic.
The canonical symplectic forms on the domain and codomain have the same formula,
\[
\Omega\left(  g,\nu,a,\delta\right)  \left(  \left(  \dot{g}_{1},\dot{\nu
}_{1},\dot{a}_{1},\dot{\delta}_{1}\right)  ,\left(  \dot{g}_{2},\dot{\nu}_{2},\dot{a}_{2},\dot{\delta}_{2}\right)  \right)  =\left\langle \dot{g}_{1},\dot{\nu}_{2}\right\rangle -\left\langle \dot{g}_{2},\dot{\nu}_{1}\right\rangle +\left\langle \nu,\left[  \dot{g}_{1},\dot{g}_{2}\right]
\right\rangle +
\left\langle \dot{a}_{1},\dot{\delta}_{2}\right\rangle
-\left\langle \dot{a}_{2},\dot{\delta}_{1}\right\rangle ,
\]
In calculating
\[
\Omega\left(  \sigma\left(  g,\nu,a,\delta\right)  \right)  \left(
\left(  \dot{g}_{1},\dot{\nu}_{1},\dot{a}_{1}+\dot{c}_{1},\dot{\delta}_{1}\right)  ,\left(  \dot{g}_{2},\dot{\nu}_{2},\dot{a}_{2}+\dot{c}_{2},\dot{\delta}_{2}\right)  \right)  ,
\]
with \(\dot{c}_{1},\dot{c}_{2}\in B^\perp\)
and \(\dot{\delta}_{1},\dot{\delta}_{2}\in B^*\) the only part containing the
\(\dot{c}_{i}\)s is \(\left\langle \dot{c}_{1},\dot{\delta}_{2}\right\rangle
-\left\langle \dot{c}_{2},\dot{\delta}_{1}\right\rangle ,\) which equals zero.
This shows that $\sigma$ is presymplectic.

It is clear that \(\sigma\) is \(G\)-equivariant. The
$H$-equivariance of \(\sigma\) follows from the $H$-equivariance of $\Gamma_\delta^*$
and the $H$-invariance of $\mu$ and $\alpha.$
%(see Remark \ref{H}). 
%and of the shift map $\overline{\Sigma}$ (see Lemma \ref{shift}).

Since $\sigma$ is $H$-equivariant, the map $\bar\sigma$ (defined above) is well-defined.
It is clearly  $G$-equivariant. Now, $K\cdot \left(\mathrm{Im}\,  \sigma\right)=
\left(G\times\mathfrak{g}^{\ast}\times A\times\left(\alpha+U\right)\right)\cap J_{K}^{-1}\left(0\right),$
which is an open subset of $J_K^{-1}(0).$ Since
$\mathrm{Im} \, \bar\sigma = \varphi \left(\mathrm{Im} \, \sigma \right)
= \varphi \left(K\cdot \mathrm{Im} \, \sigma \right),$ this implies that $\mathrm{Im} \, \bar\sigma$ is open.
Hence $\bar\sigma$ is a
 surjective submersion onto an open subset of $T^*\left(G\times_K A\right).$
For injectivity, suppose $\bar\sigma\left(\psi(w_1)\right)=\bar\sigma\left(\psi(w_2)\right),$
which is  equivalent to $\varphi\left(\sigma\left(w_1\right)\right)
=\varphi\left(\sigma\left(w_2\right)\right).$
By definition of $\varphi, $ this implies that $\sigma\left(w_1\right)=k\cdot\sigma\left(w_2\right)$
for some $k\in K.$ 
If the $A^*$ coordinates of $w_1$ and $w_2$ are $\delta_1$ and $\delta_2,$
this implies that $\alpha+\delta_1 = k\cdot \left(\alpha + \delta_2\right).$
Recall that $U$ was chosen so that the map 
$t:K\times_H U\to A^*,\,  [k,\delta]_H\mapsto k\cdot(\alpha + \delta),$
is injective.
Thus $[e,\delta_1]_H=[k,\delta_2]_H,$ which implies $k\in H.$
The $H$-equivariance of $\sigma$ implies that
$\sigma\left(w_1\right)=\sigma\left(k\cdot w_2\right),$
which implies $w_1=k\cdot w_2,$ since we have shown that $\sigma$ is injectiive. 
Thus $\psi(w_1)=\psi(k\cdot w_2)=\psi(w_2), $ which proves
injectivity of $\bar\sigma.$
Therefore $\bar\sigma$ is a bijective submersion, and hence an embedding,
onto an open subset of $T^*\left(G\times_K A\right).$
It is symplectic since $\sigma$ is presymplectic;
in fact, this is an application of Lemma \ref{fbar} at each base point.
\end{proof}

\begin{remark}
%\marginpar{details needed}
The reason for the notation \(\Gamma_{\delta}^{\ast}\)
is the following: if \(H\) is normal in \(K\) then there is a free action of \(K/H\)
on \(K\cdot\left(  \alpha+U\right)  \subset A^{\ast}.\) The Riemannian metric
defines a connection \(1\)-form \(T\left(  K\cdot\left(  \alpha+U\right)
\right)  \longrightarrow\left(  \mathfrak{k}/\mathfrak{h}\right)  \cong\left(
\mathfrak{m}\cap\mathfrak{k}\right)  \) on the principal bundle \(K\cdot\left(
\alpha+U\right)  \rightarrow K\cdot\left(  \alpha+U\right)  /\left(
K/H\right)  ,\) defined by orthogonal projection onto the vertical fibre
followed by the inverse of the infinitesimal generator map. We re-package this
connection \(1\)-form as a map \(K\cdot\left(  \alpha+U\right)  \longrightarrow
L\left(  A^{\ast},\mathfrak{m}\cap\mathfrak{k}\right)  \) and compose with the
shift map \(\left(  k,\delta\right)  \mapsto k\cdot\left(  \alpha
+\delta\right)  ,\) giving the map
\[
\Gamma:K\times U\longrightarrow L\left(  A^{\ast},\mathfrak{m}\cap
\mathfrak{k}\right)
\]
defined by
\[
\Gamma\left(  k,\delta\right)  \left(  k\cdot\left(  \xi\cdot\left(
\alpha+\delta\right)  +\varepsilon\right)  \right)  \ =\xi
\]
for every \(\xi\in\mathfrak{m}\cap\mathfrak{k}\) and \(\varepsilon\in\left(
\mathfrak{k}\cdot\alpha\right)  ^{\bot}.\) Define \(\Gamma^{\ast}:K\times
U\longrightarrow L\left(  \left(  \mathfrak{m}\cap\mathfrak{k}\right)  ^{\ast
},A\right)  \) by \(\Gamma^{\ast}\left(  k,\delta\right)  =\left(  \Gamma\left(
k,\delta\right)  \right)  ^{\ast}.\) Then for every \(\delta,\) the map
\(\Gamma^{\ast}\left(  e,\delta\right)  \) equals \(\Gamma_{\delta}^{\ast}\) as
defined in the above lemma. The proof that \(\Gamma^{\ast}\left(
e,\delta\right)  \in L\left(  \left(  \mathfrak{m}\cap\mathfrak{k}\right)  ^{\ast
},B^{\perp}\right)  \) and not just \(L\left(  \left(  \mathfrak{m}\cap
\mathfrak{k}\right)  ^{\ast},A\right)  \) is identical to the proof, in the above
lemma, that \(\Gamma_{\delta}^{\ast}\) is well-defined.
\end{remark}

\begin{remark}
\label{bijimm}The fact that every bijective immersion is a diffeomorphism
(used in the proof of the above lemma) is well known (see \cite{AMR88});
however the following short proof for finite-dimensional manifolds
seems not to be. Let \(f:M\rightarrow N\) be
a bijective immersion, and let \(m\) and \(n\) be the dimensions of \(M\) and \(N\)
respectively. Since \(f\) is an immersion, we have \(m\leq n.\) If \(m\) were
strictly less than \(n\) then every point in \(M\) would be a critical point,
which would imply (since \(f\) is surjective) that every point in \(N\) was a
critical value, contradicting Sard's theorem. Hence \(m=n,\) so \(f\) is a local
diffeomorphism at every point. Since \(f\) is bijection, it is a diffeomorphism.
\end{remark}

The composition $\tau=\bar\sigma\circ\bar\psi\circ L$ 
of the map $\bar\sigma$ from Lemma \ref{sigma} with
$\bar\psi\circ L$ from Equation \ref{E:Lpsi} is the $G$-equivariant embedding
\begin{align} \label{E:tau}
\tau
:G\times_{H}\left(  \mathfrak{m}^{\ast}\times B\times U\right)   &  \longrightarrow T^{\ast}(G\times_K A)\\
\left[  g,\nu,a,\delta\right]  _{H}  &  \longmapsto 
\varphi\left(  g,\mu+\nu
+a\diamond_{\mathfrak{h}}\delta,a-\Gamma_{\delta}^{\ast}\left(  \left.
\nu\right|  _{\mathfrak{m}\cap\mathfrak{k}}-a\diamond_{\mathfrak{m}\cap\mathfrak{k}}\delta\right)  ,\alpha+\delta\right)  , \notag
\end{align}
where $U$ and $\Gamma_\delta^*$ are as in Lemma \ref{sigma}.
%An equivalent formula is
%\begin{equation} \label{E:tau2}
%\tau:\left[  g,\nu,a,\delta\right]  _{H}  \longmapsto 
%\varphi\circ\sigma\circ l \left(  g,\nu,a,\delta\right)  , 
%\end{equation}
%where $l$ is the map from  Equation \ref{E:l}, which in the present setting has the formula 
%$\left(  g,\nu,a,\delta\right) \mapsto \left(  g,\nu + a\diamond_\mathfrak{h} \delta,a,\delta\right).$
Since $\tau$ maps $\left[e,0,0,0\right]_H$ to $\varphi(e,\mu,0,\alpha),$ 
and its image is an open subset of $T^*(G\times_K A),$ it 
is a symplectic tube.

Recall that there is a $G$-equivariant symplectomorphism $T^*s^{-1}:T^{\ast}(G\times_K V)\to T^*Q,$
for some neighbourhood $V$ of $0$ in $A.$ 
The composition of $\tau$ with $T^*s^{-1}$ will give our final result.
Unfortunately, the preimage $\tau^{-1}\left( T^{\ast}(G\times_K V)\right)$ 
doesn't have a simple description in general,
so we can only say that $T^*s^{-1}\circ\tau$ is defined on some neighbourhood of $[e,0,0,0]_H.$
However, in the special case $H=K,$ the $\Gamma_\delta^*$ term disappears, so 
the domain of $T^*s^{-1}\circ\tau$ is
%$\tau^{-1} \left(T^{\ast}(G\times_K V\right)=
$G\times_H\left(\mathfrak{m}^*\times (B\cap V) \times U\right).$
A second special case occurs if the domain of $s$ is the entire space $(G\times_K A),$
which occurs, for example, if
$K=G$ and $G$ acts linearly on $Q.$
In this case the domain of $T^*s^{-1}\circ\tau$ is simply
$G\times_H\left(\mathfrak{m}^*\times B\times U\right).$ We have proven the following:

\begin{theorem}
[Cotangent Bundle Slice Theorem]\label{mainGmuG}Let \(G\) be a Lie group acting
properly on a manifold \(Q\) and by cotangent lifts on \(T^{\ast}Q,\) which we
give the canonical cotangent symplectic form. Let \(J\) be the momentum map for
the \(G\) action, and let \(z\in T_{q}^{\ast}Q\) and \(\mu=J\left(  z\right)  .\)
Assume that \(G_{\mu}=G.\) Let \(H=G_{q}\) and \(K=G_{z},\) and let \(\mathfrak{h}\) and
\(\mathfrak{k}\) be their Lie algebras. Choose an \(H\)-invariant metric on
\(\mathfrak{g}\) and let \(\mathfrak{m}\) be the orthogonal complement to
\(\mathfrak{h.}\)  Choose a \(K\)-invariant metric on \(Q,\) and let \(A=\left(\mathfrak{g}\cdot q\right)^\perp.\)
By Palais' slice theorem, there
exists a \(K\)-invariant neighbourhood \(V\) of \(0\) in \(A\) such that
the map $s:G\times_{K}V   \rightarrow Q, \, \left[  g,a\right]  _{K}   \mapsto g\cdot\exp_{q}a,$
is a \(G\)-equivariant diffeomorphism onto a neighbourhood of $q.$
Let \(J_{K}\) be the momentum map for the cotangent lift of the twist
action of \(K\) on \(G\times A,\) and let
$\varphi:\left(  J_{K}^{-1}\left(  0\right)  \subset T^{\ast}\left(  G\times
A\right)  \right)   \rightarrow T^{\ast}\left(  G\times_{K}A\right)$
be the cotangent bundle reduction map, defined in Theorem \ref{RPCBRZ}.
Let $\alpha = \left.z\right|_A$ and \(B=\left(  \mathfrak{k}\cdot\alpha\right)^{\circ}.\) 
There exists an \(H\)-invariant neighbourhood $N$ of $(0,0,0)$ in
$\mathfrak{m}^{\ast}\oplus B\oplus B^{\ast}$
such that the map
\begin{align*}
T^{\ast}s ^{-1}\circ\bar\sigma\circ\bar\psi\circ L :G\times_{H}N
  &  \longrightarrow T^{\ast}Q\\
\left[  g,\nu,a,\delta\right]  _{H}  &  \longmapsto 
T^{\ast}s ^{-1}\circ\varphi\left(  g,\mu+\nu
+a\diamond_{\mathfrak{h}}\delta,a-\Gamma_{\delta}^{\ast}\left(  \left.
\nu\right|  _{\mathfrak{m}\cap\mathfrak{k}}-a\diamond_{\mathfrak{m}\cap\mathfrak{k}}\delta\right)  ,\alpha+\delta\right)  ,
\end{align*}
with \(\Gamma_{\delta}^{\ast}\) as in Lemma \ref{sigma}, 
is a symplectic tube around \(z.\)

If $H=K$ or $V=A$ then $N$ may be taken to equal
$\mathfrak{m}^*\times (B\cap V) \times U,$
where $U\subset B^*$ is chosen as in Lemma \ref{sigma}.
\end{theorem}

\begin{remark} \label{whynew}
There are three new aspects of this result, when compared with the general Hamiltonian slice theorem.
First, the symplectic tube is explicitly constructed, up to the cotangent lift of a Riemannian exponential
on the configuration space.
Second, we have used the cotangent-bundle-specific splitting $N_s\cong T^*B$ in the model space. 
Third, the tube has a uniqueness property; see
Lemma \ref{sigma}.
\end{remark}

\begin{remark}
\label{whyGmuG}This result depends crucially on the condition \(G_{\mu}=G,\)
for the following reasons. The isomorphism \(N_{s}\cong T^{\ast}B\)
depends on \(G_{\mu}=G\) (see Corollary \ref{NsKsubGmu}) and the isomorphism
\(G\times\mathfrak{g}_{\mu}^{\ast}\times N_{s}\cong T^{\ast}\left(  G\times
B\right)  \) depends on \(N_{s}\cong T^{\ast}B\) and also requires \(\mathfrak{g}_{\mu}^{\ast}=\mathfrak{g}^{\ast}\). The condition \(G_{\mu}=G\) is used twice in
the construction of \(\sigma:\) \ in the splitting \(\mathfrak{k}=\mathfrak{m}\)\(\cap\mathfrak{k}\oplus\mathfrak{h;}\) and in the application of
Palais' slice theorem
to the \(K\) action on \(A^{\ast},\) where it is required that \(K_{\alpha}=H\) .
Finally, $G_{\mu}=G$ is needed to guarantee that the map \(\sigma\) is symplectic, since this map
involves a shift by \(\mu\) (see the last paragraph of the proof of Lemma \ref{sigma}).
\end{remark}

When computing the symplectic tube in the cotangent bundle slice theorem 
in an example, it is easiest  to compute
the composition $T^*s^{-1} \circ \varphi$ directly, using the formula
\[
\left<T^*s^{-1} \circ \varphi (g,\nu,a,\beta), T(s\circ \pi_K) (g,\xi,a,\dot a)\right>
= \left<\nu,\xi\right> + \left<\beta, \dot a\right>,
\]
which follows directly from the definitions of the cotangent lift and the map $\varphi.$
Since the kernel of $T(s\circ \pi_K)$ is $\mathfrak{k}\cdot (G\times A),$ all elements of 
$TQ$ can be written as $T(s\circ \pi_K) (g,\xi_\perp,a,\dot a)$ for some 
$\xi_\perp\in\mathfrak{k}^\perp.$ Note that, when $\xi_\perp\in\mathfrak{k}^\perp,$ the 
$\mathfrak{k}^*$ component of $\nu$
is irrelevant in the above equation, and in particular, the term $a\diamond_\mathfrak{h} \delta$
in the formula in the cotangent bundle slice theorem is irrelevant.

A particularly simple case is when $G$ acts linearly in a vector space $Q$ and
$K=G.$ In this case, $A=T_q Q\cong Q,$ and all elements of $TQ$ can be written as 
$T(s\circ \pi_K) (g,0,a,\dot a).$  
Recalling that for linear actions,  
$s\circ \pi_K(g,a)=g \cdot (q+a),$ and identifying $A$ with $Q$, we have 
$T(s\circ \pi_K) (g,0,a,\dot a)= \left(g\cdot (q+a), g \cdot \dot a\right).$  
So the above
equation becomes
\[
\left<T^*s^{-1} \circ \varphi (g,\nu,a,\beta), \left(g\cdot (q+a), g \cdot \dot a\right)\right>
= \left<\beta, \dot a\right>,
\]
for all $\dot a \in Q,$
which is equivalent to  
\begin{align} \label{E:linKG}
T^{\ast}s^{-1}\circ\varphi\left(
g,\nu,a,\beta\right)  =\left(  g\cdot\left(  q+a\right)  ,g\cdot
\beta\right) .
\end{align}

\textbf{An alternative construction}
We now give an alternative formulation and proof of 
Theorem \ref{mainGmuG}. The new construction is more elegant but less concrete.
We will produce another $G$-equivariant
local symplectomorphism from $G\times_{H}\left(\mathfrak{m}^{\ast}\times N_{s}\right)$
to $T^{\ast}Q$ taking $\left[e,0,0\right]_{H}$ to $z,$ and then
show that it is the same as the one in Theorem \ref{mainGmuG}. 

We retain all of the definitions from earlier in this section, as well as the assumption $G_\mu=G.$
We have seen that
$G\times_{H}\left(\mathfrak{m}^{\ast}\times N_{s}\right)$
is isomorphic to $T^{\ast}\left(G\times_{H}B\right),$ so that
it suffices to find a $G$-equivariant local symplectomorphism from
$T^{\ast}\left(G\times_{H}B\right)$ to $T^{\ast}\left(G\times_{K}A\right)$
taking $\psi\left(e,0,0,0\right)$ to $\varphi\left(e,\mu,0,\alpha\right).$
It is natural to consider the cotangent lift of
some $G$-equivariant diffeomorphism from $G\times_{H}B$ to $G\times_{K}A,$
since cotangent lifts are automatically symplectic. However, the cotangent lift
of any map from $G\times_{H}B$ to $G\times_{K}A$ must map $\psi\left(e,0,0,0\right),$
which is in the zero section of $T^{\ast}\left(G\times_{H}B\right),$
to some element of the zero section of $T^{\ast}\left(G\times_{K}A\right),$
i.e., an element of the form $\varphi\left(g,0,a,0\right);$ but the
target point $\varphi\left(e,\mu,0,\alpha\right)$ is in general not of this form.
%The problem is that the original point $z\in T^{\ast}Q$
%is not, in general, in the zero section of $T^{\ast}Q,$ but in the
%{}``natural'' model of $G\times_{H}\left(\mathfrak{m}^{\ast}\times N_{s}\right)$
%as a cotangent bundle, the point $\left[e,0,0\right]_{H}$ \emph{is}
%in the zero section. 
We might try a momentum shift, but note that the shift 
$(g,\nu,a,\delta)\mapsto (g,\mu+\nu,a,\alpha+\delta)$ need not preserve $J_K^{-1}(0)$
(see Lemma \ref{imsigma}), so the ``map''
$\varphi(g,\nu,a,\delta)\mapsto \varphi(g,\mu+\nu,a,\alpha+\delta)$ is ill-defined.
%Also, it's not clear how to find a diffeomorphism
%from $G\times_H B$ to $G\times_K A.$

The idea of using cotangent lifts can be made to work, by 
``switching the roles of $A$ and $A^{\ast}$'':
modelling
$G\times_{H}\left(\mathfrak{m}^{\ast}\times N_{s}\right)$ 
as $T^{\ast}\left(G\times_{H}B^{\ast}\right)$ instead of $T^{\ast}\left(G\times_{H}B\right),$
and $T^{\ast}Q$ as $T^{\ast}\left(G\times_{K}A^{\ast}\right)$
instead of $T^{\ast}\left(G\times_{K}A\right).$
The advantages of this approach will be: (i) $z\in T^*Q$ will correspond to a point
in the zero section of $T^{\ast}\left(G\times_{K}A^{\ast}\right);$
and (ii) there is a simple local diffeomorphism from $G\times_{H}B^{\ast}$ to $G\times_{K}A^{\ast},$
namely $[g,\delta]_H\to [g,\alpha+\delta]_K$ (see Lemma \ref{F}.)

%We need to 
%find a $G$-equivariant local symplectomorphism \begin{align*}
%\tau_{2}:T^{\ast}\left(G\times_{H}B^{\ast}\right) & \longrightarrow T^{\ast}\left(G\times_{K}A^{\ast}\right)\\
%\psi_{\ast}\left(e,0,0,0\right) & \longmapsto\varphi_{\ast}\left(e,0,\alpha,0\right),\end{align*}
% where $\psi_{\ast}$ and $\varphi_{\ast}$ are the maps that appear
%in cotangent bundle reduction. It will turn out that the cotangent
%lift of the obvious map $\left[g,\delta\right]_{H}\rightarrow\left[g,\alpha+\delta\right]_{K}$
%works. 

Our starting point is the isomorphism in the following lemma, which is easily 
verified.

\begin{lemma}
\label{exchange}Let \(G\) act linearly on a vector space \(W\) and by cotangent
lifts on \(T^{\ast}W.\) With respect to the inverse dual action of \(G\) on
\(W^{\ast}\) 
and the corresponding cotangent lifted
action on \(T^{\ast}W^{\ast},\) the map
\begin{align*}
\chi:T^{\ast}W\cong W\oplus W^{\ast}  &  \longrightarrow W^{\ast}\oplus W\cong
T^{\ast}W^{\ast}\\
\left(  a,\alpha\right)   &  \longmapsto\left(  \alpha,-a\right)
\end{align*}
is a \(G\)-equivariant symplectomorphism, with respect to the standard
symplectic forms. If \(J\) and \(J_{\ast}\) are the standard momentum maps for the
\(G\) actions on \(T^{\ast}W\) and \(T^{\ast}W^{\ast}\) respectively, then \(J_{\ast
}\circ\chi=J,\) and in particular, \(J_{\ast}^{-1}\left(  0\right)  =\chi\left(
J^{-1}\left(  0\right)  \right)  .\)
\end{lemma}

It follows that\begin{align*}
\chi_{0}:T^{\ast}\left(G\times A\right)\cong G\times\mathfrak{g}^{\ast}\times A\times A^{\ast} & \longrightarrow G\times\mathfrak{g}^{\ast}\times A^{\ast}\times A\cong T^{\ast}\left(G\times A^{\ast}\right)\\
\left(g,\nu,a,\alpha\right) & \longmapsto\left(g,\nu,\alpha,-a\right)\end{align*}
 is symplectic with respect to the canonical symplectic forms, and
that $\chi_{0}\left(J_{K}^{-1}\left(0\right)\right)=J_{K,\ast}^{-1}\left(0\right),$
where $J_{K,\ast}$ is the momentum map of the cotangent lift of the
twist action of $K$ on $G\times A^{\ast}.$ Also, $\chi_{0}$
is clearly $G$-equivariant. Applying point cotangent
bundle reduction to both sides, $\chi_{0}$ induces a $G$-equivariant
symplectomorphism \[
\bar{\chi}_{0}:T^{\ast}\left(G\times_{K}A\right)\longrightarrow T^{\ast}\left(G\times_{K}A^{\ast}\right).\]
 By similar reasoning, the symplectic isomorphism\begin{align*}
\chi_{Z}:T^{\ast}\left(G\times B\right)\cong G\times\mathfrak{g}^{\ast}\times B\times B^{\ast} & \longrightarrow G\times\mathfrak{g}^{\ast}\times B^{\ast}\times B\cong T^{\ast}\left(G\times B^{\ast}\right)\\
\left(g,\nu,b,\beta\right) & \longmapsto\left(g,\nu,\beta,-b\right)\end{align*}
 maps $J_{H}^{-1}\left(0\right)$ to $J_{H,\ast}^{-1}\left(0\right),$
where $J_{H,\ast}$ is the momentum map for the cotangent-lift of
the twist action of $H$ on $G\times B^{\ast},$ and induces a $G$-equivariant
symplectomorphism \[
\bar{\chi}_{0}:T^{\ast}\left(G\times_{H}B\right)\longrightarrow T^{\ast}\left(G\times_{H}B^{\ast}\right).\]

Thus, in order to find a $G$-equivariant local diffeomorphism from
of $T^{\ast}\left(G\times_{H}B\right)$ to $T^{\ast}\left(G\times_{K}A\right)$
that maps $\psi\left(e,0,0,0\right)$ to $\varphi\left(e,0,0,\alpha\right),$
it suffices to find one, call it $\tau_{2},$ from $T^{\ast}\left(G\times_{H}B^{\ast}\right)$
to $T^{\ast}\left(G\times_{K}A^{\ast}\right)$ that maps $\psi_{\ast}\left(e,0,0,0\right)$
to $\varphi_{\ast}\left(e,0,\alpha,0\right),$ where $\psi_{\ast}$
and $\varphi_{\ast}$ are the maps that appear in cotangent bundle
reduction (Theorem \ref{RPCBRZ}), with domain and range as in the following summary diagram,
\[
\begin{tabular}
[c]{lllllll}%
$J_{H}^{-1}\left(  0\right)  $ & $\overset{\chi_{Z}}{\longrightarrow}$ &
$J_{H,\ast}^{-1}\left(  0\right)  $ & $\overset{}{\longrightarrow}$ &
$J_{K,\ast}^{-1}\left(  0\right)  $ & $\overset{\chi_{0}^{-1}}{\longrightarrow
}$ & $J_{K}^{-1}\left(  0\right)  $\\
$\downarrow\psi$ &  & $\downarrow\psi_{\ast}$ &  & $\downarrow\varphi_{\ast}$%
&  & $\downarrow\varphi$\\
$T^{\ast}\left(  G\times_{H}B\right)  $ & $\overset{\bar{\chi}_{Z}%
}{\longrightarrow}$ & $T^{\ast}\left(  G\times_{H}B^{\ast}\right)  $ &
$\overset{\tau_{2}}{\longrightarrow}$ & $T^{\ast}\left(  G\times_{K}A^{\ast
}\right)  $ & $\overset{\bar{\chi}_{0}^{-1}}{\longrightarrow}$ & $T^{\ast
}\left(  G\times_{K}A\right)  .$%
\end{tabular}
\]

The map $\tau_{2}$ will be the cotangent lift of the diffeomorphism
in the following lemma. 

%Recall that we have fixed a $K$-invariant metric on $Q,$ and we have
%used it to identify $B^{\ast}$ with $\left(\mathfrak{k}\cdot\alpha\right)^{\bot}\subset A^{\ast}.$

\begin{lemma}\label{F}
There exists an \(H\)-invariant neighbourhood \(U\) of \(\alpha\) in \(B^{\ast
}=\left(  \mathfrak{k}\cdot\alpha\right)  ^{\bot}\) such that the map
\begin{align*}
F:G\times_{H}U  &  \longrightarrow G\times_{K}\left(  K\cdot\left(
\alpha+U\right)  \right) \subset G\times_K A^* \\
\left[  g,\delta\right]  _{H}  &  \longmapsto\left[  g,\alpha+\delta\right]
_{K}\end{align*}
is a \(G\)-equivariant diffeomorphism of \(G\)-invariant neighbourhoods of
\(\left[  e,0\right]  _{H}\) and \(\left[  e,\alpha\right]  _{K}.\)
\end{lemma}

\begin{proof}
By Palais' slice theorem for linear actions (Theorem \ref{linslice}),
there exists an \(H\)-invariant neighbourhood \(U\) of \(0\) in
\(\left(  \mathfrak{k}\cdot\alpha\right)  ^{\bot}\) such that the tube
\begin{align*}
K\times_{H}U  &  \longrightarrow K\cdot\left(\alpha+U\right)  \subset A^{\ast}\\
\left[  k,\delta\right]  _{H}  &  \longmapsto k\cdot\left(  \alpha
+\delta\right)
\end{align*}
is a $K$-invariant diffeomorphism.
It follows that the map
\begin{align*}
G\times_{K} \left(K\times_{H}U  \right)
&\longrightarrow G\times_{K}\left(  K\cdot\left(\alpha+U\right)  \right) \\
\left[  g,\left[  k,\delta\right]  _H\right]  _K  &  \longmapsto\left[  g,k\cdot\left(  \alpha
+\delta\right)\right]
_{K}
\end{align*}
is a $G$-equivariant diffeomorphism.
It thus suffices to show that the following map is a $G$-equivariant diffeomorphism,
\begin{align*}
G\times_H U
&\longrightarrow G\times_{K} \left(K\times_{H}U  \right) \\
\left[  g, \delta\right]  _H  &  \longmapsto
\left[  g,\left[  e,\delta\right]  _H\right]  _K
\end{align*}
This is not hard to verify; a proof appears in \cite{OR04}.
\end{proof}

Note that the definition of $U$ is the same as in Lemma \ref{sigma}.

%We begin by assuming $\mu=0$  (we will later do a momentum shift to cover the more general
%case $G_{\mu}=G.$) 

Let $F$ be as in the previous lemma. Its cotangent lift is
the $G$-equivariant symplectomorphism \[
T^{\ast}F^{-1}:T^{\ast}\left(G\times_{H}U\right)\longrightarrow T^{\ast}\left(G\times_{K}\left(K\cdot\left(\alpha+U\right)\right)\right).\]
 Since $F$ maps $\left[e,0\right]_{H}$ to $\left[e,\alpha\right]_{K},$
it follows from the definitions of $\psi_{\ast}$ and $\varphi_{\ast}$
that $T^{\ast}F^{-1}$ 
maps $\psi_{\ast}\left(e,0,0,0\right)$
to $\varphi_{\ast}\left(e,0,\alpha,0\right).$ 
The composition $\bar\chi_0^{-1}\circ T^{\ast}F^{-1}\circ \bar\chi_Z$
maps $\psi_{\ast}\left(e,0,0,0\right)$
to $\varphi_{\ast}\left(e,0,0,\alpha\right).$ 
We compose this with the shift map
\begin{align*}
\overline{\Sigma}_{(\mu,0)}
:T^{\ast}\left(G\times_{K}A\right) & \longrightarrow T^{\ast}\left(G\times_{K}A\right)\\
\psi\left(g,\nu,a,\delta\right) & \longmapsto\psi\left(g,\mu+\nu,a,\delta\right)\end{align*}
which is easily shown to be a $G$-equivariant symplectomorphism,
by an argument similar to that in Lemma \ref{shift}.
The composition $\overline{\Sigma}_{(\mu,0)} \circ \bar\chi_0^{-1}\circ T^{\ast}F^{-1}\circ \bar\chi_Z$
maps $\psi_{\ast}\left(e,0,0,0\right)$
to $\varphi_{\ast}\left(e,\mu,0,\alpha\right).$ 
Composing with $\bar\psi\circ L,$ defined in Equation \ref{E:Lpsi}, gives a map
\begin{align*}
\overline{\Sigma}_{(\mu,0)} \circ \bar\chi_0^{-1}\circ T^{\ast}F^{-1}\circ \bar\chi_Z \circ\bar\psi\circ L
:G\times_H\left(\mathfrak{m}^{\ast}\oplus B \oplus B^*\right) 
\longrightarrow T^*\left(G\times_K V\right)
\end{align*}
taking $\left[e,0,0,0\right]_H$
to $\varphi_{\ast}\left(e,0,\alpha,0\right).$  
Finally, we compose with $T^*s^{-1}:T^*\left(G\times_K V\right) \to T^*Q,$ which forces us to restrict
the domain of the composition. The result is an alternative version of the cotangent bundle slice theorem
(Theorem \ref{mainGmuG}):

\begin{theorem} \label{altGmuG}
Under the conditions of the Theorem \ref{mainGmuG},
there exists an \(H\)-invariant neighbourhood $N$ of $(0,0,0)$ in
$\mathfrak{m}^{\ast}\oplus B\oplus B^{\ast}$
such that the map
\begin{align*}
T^*s^{-1}\circ \overline{\Sigma}_{(\mu,0)} \circ \bar\chi_0^{-1}\circ T^{\ast}F^{-1}\circ \bar\chi_Z
\circ\bar\psi\circ L
:G\times_{H}N
  &  \longrightarrow T^{\ast}Q
\end{align*}
(defined above) is a symplectic tube around \(z.\)
\end{theorem}

We will now show that the symplectic tubes in Theorems \ref{mainGmuG}
and \ref{altGmuG} are the same.
To do this, it suffices to show that 
$\overline{\Sigma}_{(\mu,0)} \circ \bar\chi_0^{-1}\circ T^{\ast}F^{-1}\circ \bar\chi_Z
=\bar\sigma,$
or equivalently,
\begin{equation}\label{E:equiv}
\bar\chi_0\circ \overline{\Sigma}_{(\mu,0)}^{-1}\circ \bar\sigma \circ \bar\chi_Z^{-1}
=T^*F^{-1}\, .
\end{equation}
It is straight-forward to check that 
\begin{equation}\label{E:sigma2}
\bar\chi_0\circ \overline{\Sigma}_{(\mu,0)}^{-1}\circ \bar\sigma \circ \bar\chi_Z^{-1}
\left(\psi_*(g,\nu,\delta,a)\right)
= \varphi_{\ast}\left(  g,\nu,\alpha+\delta,a+\Gamma_{\delta
}^{\ast}\left(  \left.  \rho\right|  _{\mathfrak{m}\cap\mathfrak{k}}+a\diamond_{\mathfrak{m}\cap\mathfrak{k}}\delta\right)  \right) ,
\end{equation}
for every $\psi_*(g,\nu,\delta,a)\in T^*(G\times_H U).$

To compute $T^*F^{-1},$ let $U$ and $F$ be as above, and
define $f:G\times U\to G\times(\alpha+U)$ by $f(g,\delta)= (g,\alpha+\delta).$
It is clear that the following diagram commutes,
\[\begin{tabular}
[c]{lllll}\(G\times U\) & \(\overset{f}{\longrightarrow}\) &
\(G\times\left(  \alpha+U\right)  \) \\
%& \(\hookrightarrow\) & \(G\times A^{\ast}\)\\
\(\pi_{H}\downarrow\) &  &\(\downarrow\pi_{K}\) \\
% &  & \(\downarrow\pi_{K}\)\\
\(G\times_{H}U\) &\(\overset{F}{\longrightarrow}\) 
&  \(G\times_K \left(K\cdot \left(  \alpha+U\right)\right)  \) \\
%& \({\hookrightarrow}\) \(G\times_{K}A^{\ast}\)
\end{tabular}
\]
where $\pi_H$ and $\pi_K$ are restrictions of the canonical projections.
Since $F$ is invertible, we have $\pi_H=F^{-1}\circ \pi_K\circ f.$
The surjectivity of $F$ also implies that every 
element of \(T\left(G\times_K \left(K\cdot \left(  \alpha+U\right)\right)\right)\)
can be expressed as 
$T\left(\pi_K\circ f\right)(g,\xi,\delta,\epsilon)
=T\pi_K(g,\xi,\alpha+\delta,\epsilon)$
for some $(g,\xi,\delta,\epsilon)\in T(G\times U).$
Hence we can compute $T^*F^{-1}$ as follows:
for any $(g,\nu,\delta,a)\in T^*(G\times U)\cap J_{K,*}^{-1}(0)$ 
and any $(g,\xi,\delta,\epsilon)\in T(G\times U),$
\begin{align} \label{E:TFcalc1}
&  \left\langle T^{\ast}F^{-1} \circ\psi_{\ast}\left(
g,\nu,\delta,a\right)  ,T\pi_{K}\left(  g,\xi,\alpha+\delta,\varepsilon
\right)  \right\rangle \\
&  =\left\langle \psi_{\ast}\left(  g,\nu,\delta,a\right)  ,T\left(
F^{-1}\circ\pi_{K}\circ f \right)  \left(  g,\xi,\delta,\varepsilon\right)
\right\rangle \notag\\
&  =\left\langle \psi_{\ast}\left(  g,\nu,\delta,a\right)  ,
T\pi_H  \left(  g,\xi,\delta,\varepsilon\right)
\right\rangle \notag\\
&  =\left\langle \nu,\xi\right\rangle +\left\langle a,\varepsilon\right\rangle
.\notag
\end{align}
Now we make the corresponding computation with the right-hand side of Equation \ref{E:sigma2},
namely
$\varphi_{\ast}\left(  g,\nu,\alpha+\delta,a
+\Gamma_{\delta}^{\ast}\left(  \left.  \rho\right|  _{\mathfrak{m}\cap\mathfrak{k}}+a\diamond_{\mathfrak{m}\cap\mathfrak{k}}\delta\right)  \right).$
Since 
$\Gamma_{\delta}^{\ast}\left(  \left.  \rho\right|  _{\mathfrak{m}\cap\mathfrak{k}}+a\diamond_{\mathfrak{m}\cap\mathfrak{k}}\delta\right)
\in B^\perp,
$
which annihilates  \(\varepsilon\in U\subset B^* ,\)
 we have,
%for any $(g,\nu,\delta,a)\in T^*(G\times U)$ and any $(g,\xi,\delta,\epsilon)\in T(G\times U),$
\begin{align} \label{E:TFcalc2}
&  \left\langle 
\varphi_{\ast}\left(  g,\nu,\alpha+\delta,a
+\Gamma_{\delta
}^{\ast}\left(  \left.  \rho\right|  _{\mathfrak{m}\cap\mathfrak{k}}+a\diamond_{\mathfrak{m}\cap\mathfrak{k}}\delta\right)  \right),
T\pi_{K}\left(  g,\xi,\alpha+\delta,\varepsilon
\right)  \right\rangle \\
&  =\left\langle \nu,\xi\right\rangle +\left\langle a
+
\Gamma_{\delta
}^{\ast}\left(  \left.  \rho\right|  _{\mathfrak{m}\cap\mathfrak{k}}+a\diamond_{\mathfrak{m}\cap\mathfrak{k}}\delta\right) ,\varepsilon\right\rangle \notag\\ 
&  =\left\langle \nu,\xi\right\rangle +\left\langle a,\varepsilon\right\rangle
. \notag
\end{align}
The calculations in Equations \ref{E:TFcalc1} and \ref{E:TFcalc2} prove that 
\begin{equation} \label{E:TFform}
T^*F^{-1}\left(\psi_*(g,\nu,\delta,a)\right)
= \varphi_{\ast}
\left(  g,\nu,\alpha+\delta,a+\Gamma_{\delta
}^{\ast}\left(  \left.  \rho\right|  _{\mathfrak{m}\cap\mathfrak{k}}+a\diamond_{\mathfrak{m}\cap\mathfrak{k}}\delta\right)  \right) ,
\end{equation}
which, together with Equation \ref{E:sigma2}, proves Equation \ref{E:equiv}.
Thus we have shown:

\begin{theorem} \label{altmainsame} The symplectic tubes in Theorems
\ref{mainGmuG} and \ref{altGmuG} are identical.
\end{theorem}

\begin{remark}
Using Lemma \ref{imsigma}, it can be easily checked that 
$\Gamma_{\delta
}^{\ast}\left(  \left.  \rho\right|  _{\mathfrak{m}\cap\mathfrak{k}}+a\diamond_{\mathfrak{m}\cap\mathfrak{k}}\delta\right) $
is the unique element $c\in B^\perp$ such that 
$\left(  g,\nu,\alpha+\delta,a+c\right)\in J_{K,*}^{-1}(0).$
Thus the calculations in Equations \ref{E:TFcalc1} and \ref{E:TFcalc2} 
show that
the formula in Equation \ref{E:TFform} is the unique one expressing
$T^*F^{-1}\left(\psi_*(g,\nu,\delta,a)\right)$ as the $\varphi_*$-image of
an element of $T^*_{(g,\alpha+\delta)} (G\times A).$
\end{remark}

\bigskip
\noindent
\textbf{Example ($SO\left(3\right)$ acting on $T^{\ast}\mathbf{R}^{3})$}
We conclude this section with a calculation in a simple example of the symplectic tube in the
cotangent bundle slice theorem (Theorem \ref{mainGmuG}). 
%and the symplectic
%isomorphism $T^{\ast}B\rightarrow N_{s}$ in Corollary \ref{NsKsubGmu}.
Consider $G=SO\left(3\right)$ acting on in the standard way on $Q=\mathbf{R}^{3},$
and by cotangent lifts on $T^{\ast}\mathbf{R}^{3}.$ The
momentum map is $\mu=J\left(q,p\right)=q\times p.$ This one example is
actually many, because we can vary the point $z=\left(q,p\right)\in T^{\ast}\mathbf{R}^{3}$
around which we construct a tube. In order
to apply Theorem \ref{mainGmuG} we require $G_{\mu}=G;$
in this case, the coadjoint action of $SO\left(3\right)$ is such
that this condition is satisfied only at $\mu=0.$ 
Thus $q$ and $p$ must be parallel, or at least one of them must be zero.
We will present the case $q=0$ and $p\neq0,$
and then state without details the results of similar calculations for the other cases.
We will implicitly use the Euclidean inner product in several places,
to define orthogonal complements and to identify spaces with their
duals. 
%\emph{Warning: the letters} $x,y,z$ \emph{will be used
%in terms such as {}``}$x$\emph{-axis'' and {}``}$yz$\emph{-plane''
%with the usual meaning in} $\mathbf{R}^{3};$ \emph{despite the notational conflicts.}

Assume $q=0$ and $p\neq0.$ Without loss of generality $z=\left(q,p\right)=\left(\left(0,0,0\right),\left(\lambda,0,0\right)\right)$
for some $\lambda\neq0$. 
%Clearly $\mu=0$ and $G_{\mu}=G$. 
We have
$K=G=SO\left(3\right),$ and $H$ is the circle group of rotations
around the $x$-axis. Since $G$ fixes $q,$ we have $A=\left(\mathfrak{g}\cdot q\right)^{\bot}=\mathbf{R}^{3}.$
Also, $\alpha=\left.z\right|_{A}=\left(\lambda,0,0\right)\in A^*.$ Since
$K=SO\left(3\right),$ the group orbit $K\cdot\alpha$ is the sphere
of radius $\lambda,$ so $B:=\left(\mathfrak{k}\cdot\alpha\right)^{\circ}$
is the $x$-axis (identifying $\left(\mathbf{R}^{3}\right)^{\ast}$
with $\mathbf{R}^{3});$ the space $B^{\ast}$ is also the $x$-axis. 
We make the standard
identification \begin{align*}
\mathfrak{g}=\mathrm{so}\left(3\right) & \longrightarrow\mathbf{R}^{3}\\
\left(\begin{array}{ccc}
 0& -\xi_{3} & \xi_{2}\\
\xi_{3} & 0 & -\xi_{1}\\
-\xi_{2} & \xi_{1}&0\end{array}\right) & \longmapsto\left(\xi_{1},\xi_{2},\xi_{3}\right).\end{align*}
 Since $H$ is the group of rotations around the $x$-axis, its Lie
algebra, $\mathfrak{h},$ is the $x$-axis.
Now $\mathfrak{m}=\mathfrak{h}^\perp$ so, 
using the Euclidean metric again, we can identify $\mathfrak{m}$ and $\mathfrak{m}^*$ with
the $yz$-plane.
We now calculate the map $\tau$  from Equation \ref{E:tau}, beginning with
the subset $U\subset B^*$ and the map $\Gamma_\delta^*$ in Lemma \ref{sigma}.
Let $U=\left(-\lambda,\infty\right)\times\left\{ \left(0,0\right)\right\} \subset B^{\ast};$
this is the largest neighbourhood in $B^{\ast}$ such that the tube
$t:K\times_{H}U\longrightarrow A^{\ast},\left[g,\delta\right]_{H}\longmapsto g\cdot\left(\alpha+\delta\right)$
is injective. 
Since $K=G,$
we have $\mathfrak{m}\cap\mathfrak{k=m.}$ The map $\Gamma_{\delta}^{\ast}:\left(\mathfrak{m}\cap\mathfrak{k}\right)^{\ast}\rightarrow B^{\perp}$
is defined by \[
\left\langle \Gamma_{\delta}^{\ast}\left(\nu\right),\xi\cdot\left(\alpha+\delta\right)+\varepsilon\right\rangle =\left\langle \nu,\xi\right\rangle ,\]
 for every $\xi\in\mathfrak{m}\cap\mathfrak{k},\nu\in\left(\mathfrak{m}\cap\mathfrak{k}\right)^{\ast},
 \delta\in U$
and $\varepsilon\in B^{\ast}.$
Let $\nu=\left(0,\nu_{2},\nu_{3}\right),\xi=\left(0,\xi_{2},\xi_{3}\right),\delta=\left(\delta_{1},0,0\right)$
and $\varepsilon=\left(\varepsilon_{1},0,0\right),$
and recall that $\alpha=\left(\lambda,0,0\right).$
 Then \[
\xi\cdot\left(\alpha+\delta\right)+\varepsilon=\xi\cdot\left(\lambda+\delta_{1},0,0\right)+\left(\varepsilon_{1},0,0\right)=\left(\varepsilon_{1},\xi_{3}\left(\lambda+\delta_{1}\right),-\xi_{2}\left(\lambda+\delta_{1}\right)\right),\]
 so $\left\langle \Gamma_{\delta}^{\ast}\left(\nu\right),\left(\varepsilon_{1},\xi_{3}\left(\lambda+\delta_{1}\right),-\xi_{2}\left(\lambda+\delta_{1}\right)\right)\right\rangle =\left\langle \nu,\xi\right\rangle .$
It follows that $\Gamma_{\delta}^{\ast}\left(\nu\right)=\left(0,\frac{\nu_{3}}{\lambda+\delta_{1}},-\frac{\nu_{2}}{\lambda+\delta_{1}}\right).$
Now $a\diamond_{\mathfrak{g}}\delta=J_K\left(a,\delta\right)=a\times\delta.$
Since $B$ is the $x$-axis and $U\subset B^{\ast},$ we
have $a\times\delta=0$ for all $a\in B$ and $\delta\in U.$
Putting these calculations together, 
\[
\tau\left[g,\nu,\left(\left(a_{1},0,0\right)\right),
\left(\delta_{1},0,0\right)\right]_H
=\varphi\left(g,\nu,\left(a_{1},\frac{\nu_{3}}{\lambda+\delta_{1}},-\frac{\nu_{2}}{\lambda+\delta_{1}}\right),\left(\lambda+\delta_{1},0,0\right)\right) \, .
\]
The symplectic tube in Theorem \ref{mainGmuG} is $T^*s^{-1}\circ \tau.$
Since $K=G=SO\left(2\right),$ we know from Equation \ref{E:linKG} that
$T^{\ast}s^{-1}\circ\varphi\left(g,\nu,a,\beta\right)=\left(g\cdot a,g\cdot\beta\right).$
So we obtain the following symplectic tube for the \(G\) action around \(z\)
(dropping the subscript-$1$'s):
\begin{align*}
T^{\ast}s^{-1}\circ\tau:G\times_{H}\left(  \mathfrak{m}^{\ast}\times
B\times U\right)   &  \longrightarrow T^{\ast}\mathbf{R}^{3}\\
\left[  g,\left(  \nu_{1},\nu_{2},\nu_{3}\right)  ,\left(  a,0,0\right)
,\left(  \delta,0,0\right)  \right]  _{H}  &  \longmapsto\left(  g\cdot\left(
a,\frac{\nu_{3}}{\lambda+\delta},-\frac{\nu_{2}}{\lambda+\delta}\right)
,g\cdot\left(  \lambda+\delta,0,0\right)  \right)
\end{align*}

The other nontrivial subcase of the $Q=\mathbf{R}^3,G=SO(3)$ example
occurs when \(q\neq0\) and \(p\parallel q;\)
it turns out that it makes no difference whether or not $p=0.$
In this case, $H=K\cong SO(2).$ Since $H=K,$ the map \(\Gamma_{\delta}^{\ast}\) 
is trivial. However, since \(K\) is neither \(G\) nor \(\left\{  e\right\}  \), the calculation of
\(T^{\ast}s^{-1}\circ\varphi\) is nontrivial, though not difficult. 
For brevity, we state only the final result for this case:
if \(z=\left(  \left(  \kappa,0,0\right)  ,\left(  \lambda,0,0\right)  \right)
\,\), then $\mathfrak{m}^*$ may be identified with the $yz$-plane, and $A$ and $B$ with the $x$-axis.
The map\begin{align*}
T^{\ast}s^{-1}\circ\tau:G\times_{H}\left(  \mathfrak{m}^{\ast}\times
(-\kappa,\infty) \times B^{\ast}\right)   &  \longrightarrow T^{\ast}\mathbf{R}^{3}\\
\left[  g,\left(  \nu_{1},\nu_{2},\nu_{3}\right)  ,\left(  a,0,0\right)
,\left(  \delta,0,0\right)  \right]  _{H}  &  \longmapsto\left(  g\cdot\left(
\kappa+a,0,0\right)  ,g\cdot\left(  \lambda+\delta,\frac{\nu_{3}}{\kappa
+a},-\frac{\nu_{2}}{\kappa+a}\right)  \right)
\end{align*}
is a symplectic tube for the \(G\) action around \(z.\)

The only remaining subcase is $q=p=0.$
In this case, $H=K=SO(3),A=B=\mathbf{R}^3$ and $U=B^*.$ The map
$T^*s^{-1}\circ \tau:G\times_G(A\times A^*)\to T^*A$ is 
the trivial one $[g,a,\delta]_G\mapsto (a,\delta).$

\section{Conclusion}

We have investigated the local structure of a cotangent bundle with
a Lie group of cotangent-lifted symmetries.
We proved a ``tangent-level'' commuting reduction result, Theorem \ref{TLCR},
and then used it in Section \ref{sectNs} to analyse the symplectic normal space.
In two special cases, we arrived at splittings of the symplectic normal space.
One of these splittings, Corollary \ref{NsKsubGmu}, applies whenever the configuration 
isotropy group $K$ is contained in the momentum isotropy group $G_\mu.$
We noted that this occurs whenever $K$ is normal in $G,$ for example when $G$ is abelian.
The splitting in Corollary \ref{NsKsubGmu} generalises the one 
given for free actions by Montgomery et al. \cite{MMR84}.
The conditions on the other splitting, in Corollary \ref{Nsalpha0}, are satisfied by
all relative equilibria of simple mechanical systems. 
In both of these special cases, the new splitting leads to a refinement of the
reconstruction equations (bundle equations), as explained at the end of Section 4.
We also noted in Section 4 two
cotangent-bundle-specific local normal forms for the symplectic reduced space,
in Theorem \ref{mainGK} and Remark \ref{CBSLNFred}.

Our main result is a cotangent bundle slice theorem, Theorem \ref{mainGmuG},
which applies at all points with fully isotropic momentum values, $G_\mu=G.$
This theorem extends the Hamiltonian slice theorem of Marle, Guillemin and
Sternberg in three ways. First, it is constructive, apart from the use of the cotangent lift of
a Riemannian exponential. Second,
it includes a cotangent-bundle-specific splitting of the symplectic normal space
(a special case of one of the first of the splittings described in the previous paragraph).
Third, our construction has a uniqueness property, contained in Lemma \ref{sigma}.
In Theorems \ref{altGmuG} and \ref{altmainsame}, we gave an alternative
contruction of the symplectic tube in the main theorem, showing that 
it is essentially a cotangent lift of a simple map between
certain twisted products.
The example presented
at the end of the Section 5 shows that our construction is feasible; we
believe that this is the first time that symplectic tubes have been
computed in an example. 

A number of open questions remain,
the most salient of which is: what happens when $\mu$ is not fully isotropic?
We have so far only been
able to formulate our cotangent bundle slice theorem for the case
of a fully-isotropic momentum value,
for reasons summarised in Remark \ref{whyGmuG}.
Our characterisation
of the symplectic normal space $N_s$ is also incomplete in the general case.
We have found a splitting of $N_s$ that applies to all relative equilibria of
simple mechanical systems, but 
what about relative equilibria of other systems?
Even for the simple mechanical case, 
what form do the reconstruction equations take if $\mu$ is not split?

One possible application of this work is to the problem
of singular cotangent bundle reduction (this was in fact our initial
motivation for this research). Local normal forms given in Section
\ref{sectNs} are a start, but do not address the global structure. 

Dynamical applications seem the most promising. To start with, the constructive
nature of the cotangent bundle slice theorem should allow us to apply theoretical results on
stability, bifurcations and persistence, such as those referred to in the Introduction,
to specific examples.
Also, the refinement of the reconstruction equations in the cotangent bundle case
may lead to extensions of the theory.
In particular, the relationship
between our splitting of the symplectic normal space at a relative equilibrium 
of a simple mechanical system,
and the splitting used in the Lagrangian Block
Diagonalisation \cite{Lew92} method for testing stability, deserves
investigation. 

\bigskip
%
%\begin{acknowledgement}
\noindent \textbf{Acknowledgment} \emph{I would like to thank Tudor Ratiu 
and Juan Pablo Ortega 
for introducing me to
this subject, and for many helpful suggestions. This article is based on part
of my Ph.D. thesis, written at the Ecole Polytechnique F\'{e}d\'{e}rale de
Lausanne under the supervision of Tudor Ratiu, supported by the Swiss National
Fund grant 21-54'138. 
Thanks also to Mark Roberts and Claudia Wulff for clarifying
the ``tangent-level commuting reduction'' viewpoint, and to M. Esmeralda Sousa-Dias,
Richard Cushman, and Miguel Rodriquez Olmos
for helpful comments on the manuscript.
%\end{acknowledgement}
}

\end{document}